# Modeling Load Redistribution Attacks in Integrated Electricity-Gas Systems

Rong-Peng Liu, *Member, IEEE,* Xiaozhe Wang, *Senior Member, IEEE,* Bo Zeng, *Member, IEEE,* and Rawad Zgheib

***Abstract*—As the share of gas-fired power generation continues to increase, it is essential to consider integrated electricity-gas systems (IEGSs). However, the development of IEGSs has led to an increased reliance on cyber infrastructures, such as communication systems, which in turn makes IEGSs more susceptible to cyberattacks. This paper investigates the paradigm and stealthy conditions of load redistribution (LR) attacks on IEGSs and proposes a bilevel mixed-integer model to identify the most severe LR attack from an economic perspective. Under a mild assumption, we prove that the proposed model does not exclude any possible upper-level attack, which differs from existing models that may yield suboptimal solutions. A modified reformulation and decomposition (R&D) algorithm is developed to solve this model in a master-subproblem framework. Particularly, we design a new subproblem to address potential infeasibility issues in the master problem. Accordingly, two types of cuts are added to the master problem for ensuring algorithm feasibility and solution optimality. Testing results validate the effectiveness of the proposed model and solution method.**

## NOMENCLATURE

### A. Sets

| | |
|---|---|
| $\mathcal{P}_d/\mathcal{G}_d$ | Set of power/gas loads. |
| $\mathcal{P}_g/\mathcal{G}_g/\mathcal{G}_w$ | Set of coal-fired/gas-fired power units/gas wells. |
| $\mathcal{P}_l/\mathcal{G}_l/\mathcal{G}_c$ | Set of power transmission lines/gas passive pipelines/gas active pipelines (compressors). |
| $\mathcal{P}_n/\mathcal{G}_n$ | Set of power/gas nodes. |
| [K] | Set $\{1, \cdots, K\}$. |

### B. Parameters

| | |
|---|---|
| $\mathrm{b}^{\mathrm{D}}_{kd}/\mathrm{b}^{\mathrm{G}}_{kg}$ | Elements in an power incidence matrix. |
| $\mathrm{c}_g/\mathrm{c}_w$ | Output costs of coal-fired power unit $g$/gas well $w$. |
| $\mathrm{c}^{\mathrm{p}}_d/\mathrm{c}^{\mathrm{g}}_d$ | Load shedding costs of power/gas load $d$. |
| $\mathrm{G}^{\min}_n/\mathrm{G}^{\max}_n$ | Minimum/maximum squared pressure bound in gas node $n$. |
| $\mathrm{G}_w$ | Maximum output capacity of gas well $w$. |
| $\mathrm{K}/\mathrm{X}^k_l$ | Number of piecewise line segments/constant piecewise endpoint. |
| $\mathrm{P}_d/\mathrm{G}_d$ | Power/gas load $d$. |
| $\mathrm{P}^{\min}_g/\mathrm{P}^{\max}_g$ | Minimum/maximum output capacity of power unit $g$. |
| $\mathrm{P}_l/\mathrm{G}_l/\mathrm{G}_c$ | Maximum transmission limit of power transmission line $l$/gas passive pipeline $l$/gas active pipeline $c$. |
| $\mathrm{W}_l$ | Weymouth constant of gas passive pipeline $l$. |
| $\alpha_c$ | Maximum compression ratio of gas active pipeline $c$. |
| $\beta_{lk}$ | Power transfer distribution factors. |
| $\tau^{\mathrm{p}}/\tau^{\mathrm{g}}$ | Maximum percentage of injected false data in the measurement of power/gas loads. |
| $\gamma_g$ | Electricity-gas conversion ratio of gas-fired power unit $g$. |

### C. Variables

| | |
|---|---|
| $s^{\mathrm{p}}_d/s^{\mathrm{g}}_d$ | Load shedding of power/gas load $d$. |
| $p_g/g_w$ | Power/gas output of power unit $g$/gas well $w$. |
| $p_l/g_l/g_c$ | Power/gas/gas flow in power transmission line $l$/gas passive pipeline $l$/gas active pipeline $c$. |
| $u_g$ | Binary unit commitment variable. 1 denotes power unit $g$ is on, and 0 otherwise. |
| $\Delta p_d/\Delta g_d$ | Injected false data in the measurement of power/ gas load $d$. |
| $\Delta p_l$ | Injected false data in the measurement of power flow in power transmission line $l$. |
| $\Delta g_l/\Delta g_c$ | Injected false data in the measurement of gas flow in gas passive pipeline $l$/gas active pipeline $c$. |
| $\Delta \pi_n$ | Injected false data in the measurement of squared pressure in gas node $n$. |
| $\pi_n$ | Squared pressure in gas node $n$. |
| $\iota^k_l/\sigma^k_l$ | Continuous/binary auxiliary variable. |

This work was supported in part by Hydro-Québec, the Institut de Valorisation des Données (IVADO), MITACS under Grant IT27493, in part by the Natural Sciences and Engineering Research Council of Canada (NSERC) under Alliance Grants ALLRP 566986-21, and in part by the Fonds de recherche du Québec – Nature et technologies (FRQNT) under Grant 334636.

R. Liu and X. Wang are with the Department of Electrical and Computer Engineering, McGill University, Montreal, QC H3A 0E9, Canada (e-mail: rongpeng.liu@mail.mcgill.ca/rpliu@eee.hku.hk; xiaozhe.wang2@mcgill.ca).

B. Zeng is with the Department of Industrial Engineering and the Department of Electrical and Computer Engineering, University of Pittsburgh, Pittsburgh, PA 15106 USA (e-mail: bzeng@pitt.edu).

R. Zgheib is with Hydro-Quebec Research Institute (IREQ) (email: zgheib.rawad2@hydroquebec.com).

## I. INTRODUCTION

THE share of gas-fired power generation has been continuing to increase over the last decade, contributing to 38.3% and 40.0% of the aggregated electricity in the U.S. [1] and U.K. [2], respectively, in 2021. It is essential to consider electricity and gas systems holistically as an integrated electricity and gas system (IEGS). However, the development of IEGSs has led to an increased reliance on cyber infrastructures, such as communication systems for enabling the coordinated operation of the power and gas subsystems in an IEGS. It is the increased reliance on cyber infrastructures that in turn makes IEGSs more susceptible to cyberattacks. Many cyberattacks on the energy sector have been recorded, e.g., the 2015 Ukraine power system hack [3] and the 2021 Colonial Pipeline ransomware cyberattack [4]. For the latter, the company paid millions of dollars (to a hacker group) to stop this cyberattack [5]. While efforts have been made to establish cybersecurity directives and protocols to mitigate cyber threats [6], it is still concerning to note that a significant portion of the natural gas subsector has a high risk of suffering ransomware attacks due to a high level of system automation [7]. Given the substantial share of electricity gener-ated by gas-fired units and the importance of IEGSs, systematic studies of cyberattacks on IEGSs are highly demanded.

Researchers have made intensive efforts to investigate various cyberattack paradigms on power systems [8]-[18]. Among them, references [8]-[10] studied false data injection attacks (FDIAs), where intruders inject elaborate false data into power system measurements (e.g., direct injection into the meter) without being detected by bad data detection (BDD), thereby causing incorrect estimation of system states. Reference [8] proposed FDIAs against direct current (DC) power system state estimation. Reference [9] extended FDIAs to eroding alternating current power system state estimation. Reference [10] designed FDIAs that require only partial system information.

As a type of FDIA, load redistribution (LR) attacks on power systems received huge attention [11]-[14]. Unlike FDIAs, LR attacks do not compromise generator output measurements, as falsified generator output measurements can be detected due to the direct communication between multiple control centers [11]. In other words, LR attacks are a practical type of FDIAs that are challenging to detect. In addition, the study of LR attacks also focused on identifying the most severe attack with a specific objective, e.g., operation cost. Compared with FDIAs, LR attacks may have a more substantial impact on operational decisions, resulting in significant consequences, e.g., economic losses. Specifically, reference [11] studied LR attacks against DC power system state estimation and proposed a bilevel model to identify the most severe attack aiming to maximize operation costs. Reference [12] proposed local LR attacks which required only partial system information. Reference [13] confirmed the potential of LR attacks for inducing severe cascading failures. Reference [14] designed stealthier LR attacks, i.e., dummy data attacks, as classical LR attacks could be detected by clustering methods. Other types of cyberattacks (on power systems), e.g., data integrity attacks [15], coordinated cyber-physical attacks [16], [17], and denial-of-service attacks [18], were also studied. These cyberattack paradigms motivate researchers to develop detection and mitigation methods [19], [20]. This work is particularly interested in LR attacks due to their practical nature and severity, especially in the identification of the most severe attack scenarios.

Although extensive efforts have been dedicated to investigating cyberattack paradigms on power systems, limited works studied IEGS cybersecurity [21]-[27]. Reference [21] proposed an adversarial model to analyze the vulnerability of IEGSs under cyberattacks that endeavor to gain access to gas compressor control systems. References [22]-[26] focused on FDIAs on IEGSs. Specifically, reference [22] designed a coordinated day-ahead scheme against FDIAs on gas system components in an IEGS. Reference [23] developed a trilevel model to mitigate the impact of FDIAs on the gas flows supplied to gas-fired units in an IEGS. Reference [24] studied detection and prevention strategies against FDIAs on power system components and gas pipeline monitor systems in an IEGS. Reference [25] adopted a two-stage distributionally robust optimization model to study mitigation strategies to reduce the impact of FDIAs on IEGSs. The same model was used in [26] to make cyber-secured operation strategies against FDIAs on integrated electricity-gas-water systems. Regarding LR attacks on IEGSs, reference [27] adopted the same model as [25] to derive preventive strategies for distribution-level IEGSs considering energy hubs against LR attacks on power system components. It is important to note that reference [27] solely studied LR attacks on power system components in IEGSs but ignored possible attacks on gas systems. Moreover, it did not consider interdependency between power and gas systems in IEGSs and may make LR attacks (on IEGSs) not stealthy. In fact, extending LR attacks from power systems to *entire IEGSs* is challenging due to (implicit) infeasibility issues and nonconvex gas transmission equations [21].

Although cyberattacks pose great threats to IEGSs [7], to the best of our knowledge, no existing works studied LR attacks on entire IEGSs, let alone identifying the *most severe attacks*. We summarize the research gaps: i) the paradigm of LR attacks on IEGSs has not been systematically studied, ii) no conditions are proposed to ensure the stealthiness of the attacks, iii) nonconvex gas transmission equations incur challenges in designing LR attacks on IEGSs and identifying the most severe attack.

To address these gaps, we conduct a comprehensive exploration of LR attacks on IEGSs, including both modeling and solution methods. Overall, our contributions are twofold:

1) This paper studies LR attacks on IEGSs. For the first time, we develop the paradigm of LR attacks on an entire IEGS and give conditions to ensure the stealthiness of the LR attacks. Furthermore, we propose a bi-level mixed-integer linear program (MILP) to identify the most severe LR attack from an economic perspective, where i) unit commitment (UC) variables are incorporated, and ii) nonconvex gas transmission constraints are piecewisely linearized by a piecewise linear (PWL) method. We prove that, under a mild assumption, the proposed model does not exclude any possible upper-level attack. This differentiates our model from existing bilevel LR attack models that may ignore possible LR attacks and yield suboptimal solutions.

2) This paper designs a modified reformulation and decomposition (R&D) algorithm, building upon the R&D algorithm [28], for solving the proposed bilevel MILP in a master-subproblem framework. Particularly, we propose a new subproblem to address potential infeasibility issues in the master problem (by checking if all possible upper-level attacks are feasible under fixed lower-level binary variables). Accordingly, two types of cuts are added to the master problem for ensuring both the algorithm feasibility and the solution optimality. We prove that the modified R&D algorithm converges within a finite number of iterations. Simulation results validate the effectiveness of the proposed model and solution method.

Note that there are many algorithms for solving bilevel linear programs [29], but the algorithms for solving bilevel MILPs are limited [30], [31]. Reference [30] proposed a value function-based algorithm for the bilevel MILPs with integer upper-level problems. Reference [31] required that the lower-level problem should not contain upper-level continuous variables. In view of the mathematical formulation of the proposed bilevel MILP (referred to as **O-M** in Section II.C), neither algorithm is applicable. Hence, the modified R&D algorithm is specifically designed to address the proposed bilevel MILP.

The remainder of this paper is structured as follows. Section

II analyzes the LR attacks on IEGSs and proposes a bilevel model for deriving the most threatening attack. Section III develops the modified R&D solution method. Section IV conducts numerical experiments. Section V concludes this paper.

## II. Problem Formulation

First, we clarify assumptions and simplifications for IEGSs.

1) This work aims at transmission-level IEGSs and uses i) an DC power flow model to depict power flow in power transmission lines [8], [11] and ii) a Weymouth equation to model gas flow in long-distance high-pressure gas passive pipelines [21].

2) Each gas active pipeline is equipped with a compressor. For depicting gas flow in gas active pipelines, we adopt the simplified compressor model to avoid (unnecessary) computational burdens [21]-[27], [32].

3) Intruders have complete information about IEGS topology and parameters and can manipulate system measurements. Although it is restrictive, this assumption is broadly adopted, e.g., [21]-[27]. For some situations, we note that intruders can still manage to obtain sufficient system data stored in control centers [8] through real-world cyberattacks [3], [4]. Certainly, future efforts should consider relaxing this assumption.

### A. Preliminary: LR Attacks on Power Systems

State estimation is widely adopted to estimate power system states by means of measurements. Intruders can conduct cyberattacks by injecting well-designed false data into measurements without being detected by BDD, causing incorrect estimations of system states. Thus, it is imperative to distinguish measurable variables, i.e., measurements, from state variables. For DC power systems, Table I presents the details [9].

TABLE I
Measurements and State Variables in DC Power Systems

| Variable | Description | Measurement | State variable |
|---|---|---|---|
| $p_n$ ($n \in \mathcal{P}_n$) | (active) power injection | ✓ | × |
| $p_l$ ($l \in \mathcal{P}_l$) | (active) power flow | ✓ | × |
| $\theta_n$ ($n \in \mathcal{P}_n$) | phase angle | * | ✓ |

*: The phase angles are measurable only for the power nodes equipped with PMUs.

An LR attack on power systems is implemented by injecting false data into *power load measurements* while generator output measurements keep unchanged [11]. Other measurements are manipulated correspondingly for keeping consistent with the falsified power load measurements and ensuring the stealthiness of the attack. With that said, LR attacks (on power systems) are stealthy attacks satisfying

$$\sum_{d\in\mathcal{P}_d} \Delta p_d = 0, \tag{1a}$$

$$-\tau^{\mathrm{p}} \mathrm{P}_d \le \Delta p_d \le \tau^{\mathrm{p}} \mathrm{P}_d, \quad \forall d \in \mathcal{P}_d, \tag{1b}$$

$$\Delta p_l = \sum_k \beta_{lk} \big(\sum_{d\in\mathcal{P}_d} \mathrm{b}^{\mathrm{D}}_{kd} \Delta p_d\big), \quad \forall l \in \mathcal{P}_l. \tag{1c}$$

Equation (1a) enforces the power supply-demand balance. Constraint (1b) indicates that the injected false data (in power load measurements) should not exceed certain thresholds. Otherwise, the (falsified) data may be suspected due to the large deviation from load forecasting values. Equation (1c) quantifies the false data injected into power flow measurements.

In general, there is more than one attack strategy satisfying (1). Most of previous works (e.g., [11]) focused on identifying the most severe LR attacks that maximize system costs by solving the bilevel optimization model (2).

$$\max_{\Delta\mathbf{p}_D,\, \mathbf{p}_G^0,\, \mathbf{s}_D^{\mathrm{p},0}} \sum_{g\in\mathcal{P}_g} \mathrm{c}_g p_g^0 + \sum_{d\in\mathcal{P}_d} \mathrm{c}^{\mathrm{p}}_d s_d^{\mathrm{p},0} \tag{2a}$$

$$\text{s.t. (1a), (1b),} \tag{2b}$$

$$(\mathbf{p}_G^0,\, \mathbf{s}_D^{\mathrm{p},0}) \in \arg\min_{\mathbf{p}_G,\, \mathbf{s}_D^{\mathrm{p}}} \sum_{g\in\mathcal{P}_g} \mathrm{c}_g p_g + \sum_{d\in\mathcal{P}_d} \mathrm{c}^{\mathrm{p}}_d s^{\mathrm{p}}_d \tag{2c}$$

$$\text{s.t. } \mathrm{P}_g^{\min} \le p_g \le \mathrm{P}_g^{\max}, \quad \forall g \in \mathcal{P}_g, \tag{2d}$$

$$-\mathrm{P}_l \le \sum_k \beta_{lk} \big(\sum_{g\in\mathcal{P}_g} \mathrm{b}^{\mathrm{G}}_{kg} p_g + \sum_{d\in\mathcal{P}_d} \mathrm{b}^{\mathrm{D}}_{kd} (\mathrm{P}_d + \Delta p_d - s^{\mathrm{p}}_d)\big) \le \mathrm{P}_l, \quad \forall l \in \mathcal{P}_l, \tag{2e}$$

$$\sum_{g\in\mathcal{P}_g} p_g = \sum_{d\in\mathcal{P}_d} (\mathrm{P}_d - s^{\mathrm{p}}_d), \tag{2f}$$

$$0 \le s^{\mathrm{p}}_d \le \mathrm{P}_d + \Delta p_d, \quad \forall d \in \mathcal{P}_d. \tag{2g}$$

$\Delta\mathbf{p}_D = col(\Delta p_d)$, $d \in \mathcal{P}_d$. Function $col(\cdot)$ denotes the mapping that reshapes scalars and/or vectors as one vector, e.g., $col(\mathbf{a}, \mathrm{b}, \mathbf{c}, \mathrm{d}) = [\mathbf{a}^{\mathrm{T}}\ \mathrm{b}\ \mathbf{c}^{\mathrm{T}}\ \mathrm{d}]^{\mathrm{T}}$. All vectors in this paper are column vectors, and the superscript "T" denotes the transpose. Hereinafter, all vector variables with upper-case subscripts are generated using the same rule without further statements. The upper-level attack model consists of (2a) and (2b). Objective (2a) maximizes the power system costs, i.e., the sum of generator output and power load shedding costs. Constraint (2b) ensures the stealthiness of the LR attack. The lower-level optimal power flow model is composed of (2c)-(2g), which aims to minimize system costs (2c). Constraint (2d) states the generator output boundaries. Constraint (2e) limits the power flow in power transmission lines. Equation (2f) is the power balance constraint. Constraint (2g) specifies the power load shedding range. After deriving the optimal $\Delta\mathbf{p}_D$ by (2), we can calculate $\Delta\mathbf{p}_L$ by (1c). Note that the upper level of (2) does not contain (1c), as falsified phase angle measurements (if applicable) can be derived by (1c) after obtaining falsified power load measurements.

The bilevel linear program (LP) (2) can be reformulated as an MILP by Karush-Kuhn-Tucker (KKT) conditions, which is solvable by off-the-shelf solvers, e.g., Gurobi. Other methods, such as a value-function-based method [29], are also applicable to solving (2).

### B. LR Attacks on IEGSs

In this subsection, we will extend LR attacks on power systems to IEGSs and propose a bilevel MILP for identifying the most severe LR attack. It will be shown that this extension is nontrivial due to (implicit) feasibility issues and nonconvex gas transmission equations. Hereinafter, unless stated otherwise, an *LR attack* refers to the LR attack on IEGSs. First, by referring to [33], we divide gas system variables into measurable variables, i.e., measurements, and state variables. Table II presents the details. Gas injection is derived by

$$g_n = \sum_{w\in\mathcal{G}_w} g_{w(n)} - \sum_{d\in\mathcal{G}_d} \mathrm{G}_{d(n)} - \sum_{g\in\mathcal{G}_g} \gamma_g p_{g(n)}, \quad \forall n \in \mathcal{G}_n,$$

where $w(n)$, $d(n)$, and $g(n)$ represent the gas wells, gas loads, and gas-fired generators connected to gas node $n$, respectively.

TABLE II
MEASUREMENTS AND STATE VARIABLES IN GAS SYSTEMS

| Variable | Description | Measurement | State variable |
|---|---|---|---|
| $g_n$ ($n \in \mathcal{G}_n$) | gas injection | ✓ | × |
| $g_l$ ($l \in \mathcal{G}_l$) | gas flow in gas passive pipeline | ✓ | × |
| $g_c$ ($c \in \mathcal{G}_c$) | gas flow in gas active pipeline | ✓ | ✓ |
| $\pi_n$ ($n \in \mathcal{G}_n$) | squared nodal pressure | ✓ | ✓ |

The LR attack on IEGSs aims to inject false data into *power and gas load measurements* and keeps power generator and gas well output measurements fixed. Accordingly, intruders manipulate the remaining measurements to ensure the stealthiness of the attack. Similar to power systems, IEGS measurements are derived by meters [34]. For achieving an LR attack, intruders may directly compromise IEGS meters or hack the devices which store measurements. Mathematically, it should satisfy

$$\text{(1a)-(1c)}, \tag{3a}$$

$$\sum_{d\in\mathcal{G}_d} \Delta g_d = 0, \tag{3b}$$

$$-\tau^{\mathrm{g}} \mathrm{G}_d \le \Delta g_d \le \tau^{\mathrm{g}} \mathrm{G}_d, \quad \forall d \in \mathcal{G}_d, \tag{3c}$$

$$\sum_{d\in\mathcal{G}_d} \Delta g_{d(n)} = \sum_{l\in\mathcal{G}_l} \Delta g_{l(n)} + \sum_{c\in\mathcal{G}_c} \Delta g_{c(n)}, \quad \forall n \in \mathcal{G}_n, \tag{3d}$$

$$(\tilde{\mathrm{g}}_l + \Delta g_l)|\tilde{\mathrm{g}}_l + \Delta g_l| - \tilde{\mathrm{g}}_l|\tilde{\mathrm{g}}_l| = \mathrm{W}_l(\Delta\pi_{m(l)} - \Delta\pi_{n(l)}), \quad \forall l \in \mathcal{G}_l. \tag{3e}$$

Constraint (3a) refers to stealthy conditions (1a)-(1c). Equation (3b) denotes the gas supply-demand balance. Constraint (3c) states that the injected false data in gas load measurements are within a range. Equations (3d) and (3e) ensure the consistency between all (falsified) measurements in gas systems. Function $|\cdot|$ denotes the absolute value mapping. Specifically, equations (3d) and (3e) are derived from the nodal gas balance relation and the Weymouth equation, respectively. (The formulation is shown later.) Based on [9], [10], the $\tilde{\mathrm{g}}_l$ in (3e) is set to the (true) measurement value of the gas flow in gas passive pipeline $l$.

Similar to the LR attack on power systems, intruders aim to design an LR attack to maximize the operation costs of an IEGS. For this purpose, the LR attack on IEGSs can be formulated as the following bilevel model (from intruders' perspective).

$$\max_{\mathbf{x},\mathbf{y}'_0} \sum_{g\in\mathcal{P}_g} \mathrm{c}_g p_g^0 + \sum_{d\in\mathcal{P}_d} \mathrm{c}_d^{\mathrm{p}} s_d^{\mathrm{p},0} + \sum_{w\in\mathcal{G}_w} \mathrm{c}_w g_w^0 + \sum_{d\in\mathcal{G}_d} \mathrm{c}_d^{\mathrm{g}} s_d^{\mathrm{g},0} \tag{4a}$$

$$\text{s.t. } \text{(1a), (1b), (3b), (3c)}, \tag{4b}$$

$$\mathbf{y}'_0 \in \arg\min_{\mathbf{y}'} \sum_{g\in\mathcal{P}_g} \mathrm{c}_g p_g + \sum_{d\in\mathcal{P}_d} \mathrm{c}_d^{\mathrm{p}} s_d^{\mathrm{p}} + \sum_{w\in\mathcal{G}_w} \mathrm{c}_w g_w + \sum_{d\in\mathcal{G}_d} \mathrm{c}_d^{\mathrm{g}} s_d^{\mathrm{g}} \tag{4c}$$

$$\text{s.t. } \mathrm{P}_g^{\min} \le p_g \le \mathrm{P}_g^{\max}, \quad \forall g \in \mathcal{P}_g \cup \mathcal{G}_g, \tag{4d}$$

$$-\mathrm{P}_l \le \sum_k \beta_{lk} \Big( \sum_{g\in\mathcal{P}_g\cup\mathcal{G}_g} \mathrm{b}_{kg}^{\mathrm{G}} p_g + \sum_{d\in\mathcal{P}_d} \mathrm{b}_{kd}^{\mathrm{D}} (\mathrm{P}_d + \Delta p_d - s_d^{\mathrm{p}}) \Big) \le \mathrm{P}_l, \quad \forall l \in \mathcal{P}_l, \tag{4e}$$

$$\sum_{g\in\mathcal{P}_g\cup\mathcal{G}_g} p_g = \sum_{d\in\mathcal{P}_d} (\mathrm{P}_d - s_d^{\mathrm{p}}), \tag{4f}$$

$$0 \le s_d^{\mathrm{p}} \le \mathrm{P}_d + \Delta p_d, \quad \forall d \in \mathcal{P}_d, \tag{4g}$$

$$0 \le g_w \le \mathrm{G}_w, \quad \forall w \in \mathcal{G}_w, \tag{4h}$$

$$\mathrm{G}_n^{\min} \le \pi_n \le \mathrm{G}_n^{\max}, \quad \forall n \in \mathcal{G}_n, \tag{4i}$$

$$g_l |g_l| = \mathrm{W}_l (\pi_{m(l)} - \pi_{n(l)}), \quad \forall l \in \mathcal{G}_l, \tag{4j}$$

$$\pi_{n(c)} \le \alpha_c \pi_{m(c)}, \quad \forall c \in \mathcal{G}_c, \tag{4k}$$

$$-\mathrm{G}_l \le g_l \le \mathrm{G}_l, \quad \forall l \in \mathcal{G}_l, \tag{4l}$$

$$0 \le g_c \le \mathrm{G}_c, \quad \forall c \in \mathcal{G}_c, \tag{4m}$$

$$\sum_{w\in\mathcal{G}_w} g_{w(n)} + \sum_{l\in\mathcal{G}_l} g_{l(n)} + \sum_{c\in\mathcal{G}_c} g_{c(n)} = \sum_{d\in\mathcal{G}_d} (\mathrm{G}_{d(n)} + \Delta g_{d(n)} - s_{d(n)}^{\mathrm{g}}) + \sum_{g\in\mathcal{G}_g} \gamma_g p_{g(n)}, \quad \forall n \in \mathcal{G}_n, \tag{4n}$$

$$0 \le s_d^{\mathrm{g}} \le \mathrm{G}_d + \Delta g_d, \quad \forall d \in \mathcal{G}_d. \tag{4o}$$

$\mathbf{x} = col(\Delta\mathbf{p}_D, \Delta\mathbf{g}_D)$ and $\mathbf{y}' = col(\mathbf{p}_G, \mathbf{s}_D^{\mathrm{p}}, \mathbf{g}_W, \mathbf{g}_L, \mathbf{g}_C, \boldsymbol{\pi}_N, \mathbf{s}_D^{\mathrm{g}})$. The upper-level model (4a), (4b) determines the intruders' action. Objective (4a) maximizes IEGS costs, i.e., the sum of power generator output, power load shedding, gas well output, and gas load shedding costs. Constraint (4b) reveals all possible LR attacks. We define a *possible LR attack* as an LR attack satisfying (4b). The lower-level model (4c)-(4o) provides the economic dispatch decision from the IEGS operators' perspective *under* the *falsified load measurements*. Objective (4c) minimizes the IEGS costs. Power system operation constraints (4d)-(4g) are almost the same as (2d)-(2g). The only difference is that (4d)-(4g) incorporate gas-fired power generators. Gas system operation constraints incorporate (4h)-(4o). Specifically, constraint (4h) exhibits gas well output capacity. Constraint (4i) states gas nodal pressure limits. Equation (4j), Weymouth equations [22], simulates the gas flow in gas passive pipelines, where nodes $m(l)$ and $n(l)$ are connected by gas passive pipeline $l$. Constraint (4k), the simplified compressor model [22], [32], regulates nodal pressures for gas active pipelines. Constraints (4l) and (4m) indicate the transmission capacity of gas passive and active pipelines, respectively. Equation (4n) is the gas supply-demand balance constraint, where $w(n)$, $l(n)$, $c(n)$, $d(n)$, and $g(n)$ denote the gas wells, gas flow in gas passive and active pipelines, gas loads, and gas-fired power generators connected to gas node $n$, respectively. Constraint (4o) is the gas load shedding range. After deriving the optimal $\Delta\mathbf{p}_D$ and $\Delta\mathbf{g}_D$ by (4), we can calculate $\Delta\mathbf{p}_L$, $\Delta\mathbf{g}_L$, $\Delta\mathbf{g}_C$, $\Delta\boldsymbol{\pi}_N$, by (1c), (3d), and (3e). That is why constraints (4b) do not consider stealthy conditions (1c), (3d), and (3e).

It is important to highlight that the bilevel model (4), a direct extension of the bilevel LR attack model (on power systems [11]), may fail to identify (possible) severe LR attacks. An illustrative example is shown in Appendix A, where a significantly severe attack scenario is disregarded by model (4) due to the infeasibility of its lower-level problem. To handle this issue, we propose an augmented bilevel model in the next subsection, which integrates UC variables. The augmented model is guaranteed to be feasible for all possible LR attacks.

### C. Augmented Bilevel Model

In fact, it is the constraint (4d) that may incur the infeasibility of the lower-level problem (4c)-(4o) under some (possible) LR attacks. To see this, we define the *original LR attack region* $\mathcal{R}_{\mathrm{OLR}}$ and the *inducible LR attack region* $\mathcal{R}_{\mathrm{ILR}}$ of (4) as

$$\mathcal{R}_{\mathrm{OLR}} = \{\mathbf{x} \in \mathbb{R}^n : \text{(1a), (1b), (3b), (3c)}\},$$

$$\mathcal{R}_{\mathrm{ILR}} = \{\mathbf{x} \in \mathbb{R}^n : \exists \mathbf{y}' \in \mathbb{R}^{q'} \text{ with } (\mathbf{x}, \mathbf{y}') \in \{\mathcal{R}_{\mathrm{OLR}} \cap \text{(4d)-(4o)}\}\}.$$

Evidently, any $\mathbf{x}$ in $\mathcal{R}_{\mathrm{OLR}}$ is a stealthy LR attack. Intruders can select any $\mathbf{x}$ in $\mathcal{R}_{\mathrm{OLR}}$ freely as an LR attack on IEGSs. On the other hand, the selected (upper-level) LR attack should be independent of lower-level decisions, i.e., $\mathcal{R}_{\text{ILR}}$ is expected to be equal to $\mathcal{R}_{\text{OLR}}$. However, according to the illustrative example in Appendix A, this may not hold for (4) due to the potential infeasibility of its lower level, which renders $\mathcal{R}_{\text{ILR}} \subseteq \mathcal{R}_{\text{OLR}}$.

One remedy method is to add UC variables to (4d), i.e.,

$$u_g \mathrm{P}_g^{\min} \le p_g \le u_g \mathrm{P}_g^{\max}, \quad \forall g \in \mathcal{P}_g \cup \mathcal{G}_g. \tag{5}$$

Constraint (5) allows generator outputs to be zero, i.e., generator shutdown. This constraint ensures that augmented power subsystem feasible region (4e)-(4g), (5) is not empty under any possible LR attacks on the power subsystem of an IEGS. The above claim is obvious, as $col(\mathbf{p}_G = \mathbf{0}, \mathbf{u}_G = \mathbf{0}, \mathbf{s}_D^{\mathrm{p}} = \mathbf{P}_\mathrm{D} + \Delta\mathbf{p}_D)$ is always a feasible solution (i.e., shutting down all generators and cutting off all power loads) for (4e)-(4g) and (5). By replacing (4d) with (5), we eliminate the lower-level infeasibility in (4) caused by the power subsystem constraints. Next, we check if there are other infeasibility issues in the lower level of (4) caused by gas subsystem constraints (4h)-(4o). In fact, under Assumption 1, the augmented lower-level problem (4c), (4e)-(4o), (5) is always feasible for any $\mathbf{x} \in \mathcal{R}_{\text{OLR}}$.

*Assumption 1*: The feasible region (4i)-(4k) is not empty when all $g_l$, $l \in \mathcal{G}_l$, in (4j) are set to zero.

Assumption 1 includes (but is not limited to) the practical operation status where an entire gas system is shut off, i.e., no gas-fired generator outputs and no loads. This status satisfies assumption 1 and can occur physically due to extreme events and/or severe faults [35], [36], e.g.. Thus, this assumption is considered mild in practice. Under this assumption, we have:

*Proposition 1*: If Assumption 1 holds, the augmented lower-level problem (4c), (4e)-(4o), (5) is always feasible for any possible (upper-level) LR attack and has a finite optimal value.

TABLE III
FEASIBLE SOLUTION

| Variable | $\mathbf{p}_G$ | $\mathbf{u}_G$ | $\mathbf{s}_D^{\mathrm{p}}$ | $\mathbf{g}_W$ | $\mathbf{g}_L$ | $\mathbf{g}_C$ | $\boldsymbol{\pi}_N$ | $\mathbf{s}_D^{\mathrm{g}}$ |
|---|---|---|---|---|---|---|---|---|
| Value | $\mathbf{0}$ | $\mathbf{0}$ | $\mathbf{P}_\mathrm{D}+\Delta\mathbf{p}_D$ | $\mathbf{0}$ | $\mathbf{0}$ | $\mathbf{0}$ | ** | $\mathbf{G}_\mathrm{D}+\Delta\mathbf{g}_D$ |

**: The values of $\pi_n$, $n \in \mathcal{G}_n$, satisfy Assumption 1.

The proof for Proposition 1 is straightforward, as the solution presented in Table III is always feasible for the lower-level problem (4c), (4e)-(4o), (5) with any possible attack $\Delta\mathbf{p}_D$ and $\Delta\mathbf{g}_D$, where $\mathbf{P}_\mathrm{D} = col(\mathrm{P}_d)$, $d \in \mathcal{P}_d$, and $\mathbf{G}_\mathrm{D} = col(\mathrm{G}_d)$, $d \in \mathcal{G}_d$. Proposition 1 states that no possible (upper-level) LR attacks are excluded due to the lower-level infeasibility.

The augmented bilevel LR attack model (4a)-(4c), (4e)-(4o), (5) is a bilevel program with a mixed-integer nonconvex lower level. Particularly, nonconvex Weymouth equations (4j) incur significant computational difficulties. We adopt an incremental PWL method [37] to piecewisely linearize them. For gas passive pipeline $l$, constraint (4j) is replaced by

$$h(\mathrm{X}_l^1) + \sum_{k \in [\mathrm{K}]} (h(\mathrm{X}_l^{k+1}) - h(\mathrm{X}_l^k)) \iota_l^k = \pi_{m(l)} - \pi_{n(l)}, \tag{6a}$$

$$g_l = \mathrm{X}_l^1 + \sum_{k \in [\mathrm{K}]} (\mathrm{X}_l^{k+1} - \mathrm{X}_l^k) \iota_l^k, \tag{6b}$$

$$\iota_l^{k+1} \le \sigma_l^k \le \iota_l^k, \;\; k \in [\mathrm{K}-1], \tag{6c}$$

$$0 \le \iota_l^k \le 1, \;\; k \in [\mathrm{K}], \tag{6d}$$

where function $h(\mathrm{X}) = \mathrm{X}|\mathrm{X}|/\mathrm{W}_l$. Although we approximately simplify the mixed-integer nonconvex Weymouth equation (4j) into simpler mixed-integer linear constraints (6), the incremental PWL method also introduces errors. Fig. 1 shows an example of piecewisely linearizing a Weymouth equation by the incremental PWL method. In this case, a Weymouth equation is separated into eight segments. Here, we always require that the origin (0, 0) is an endpoint of a segment, e.g., ($\mathrm{X}^5$, $h(\mathrm{X}^5)$) in Fig. 1. Accordingly, we analyze the errors in the following aspects.

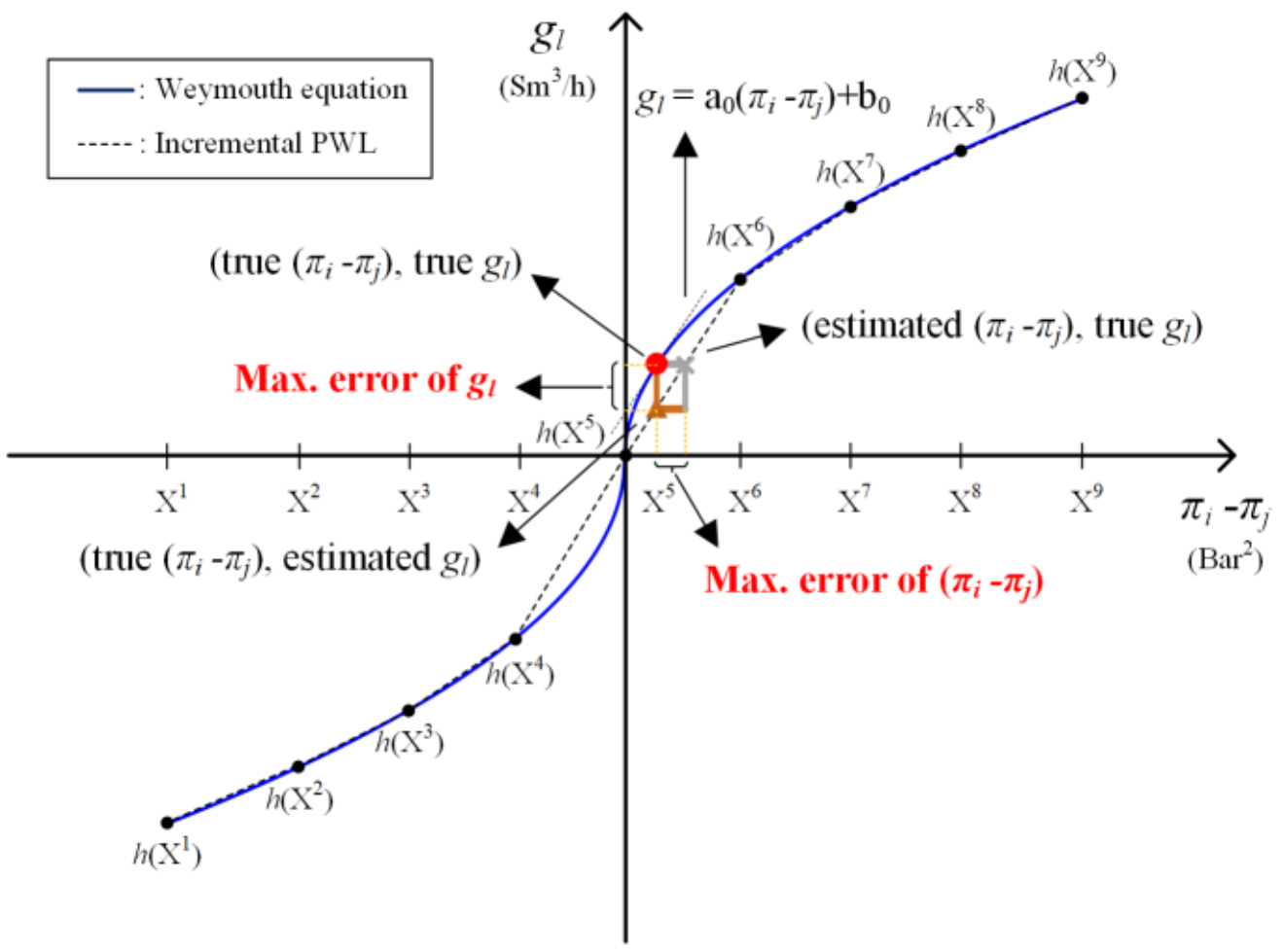

Fig. 1. Incremental PWL of a Weymouth equation.

1) Trend. According to Fig. 1, the incremental PWL method is more liable to incur smaller errors when the squared pressure difference, i.e., $(\pi_i - \pi_j)$, is large. Physically, this scenario corresponds to the IEGS with heavy gas loads, as gas pipelines usually transmit more gas, i.e., larger $g_l$ in (4j), and thus require larger squared pressure differences.

2) Range. We assume that $(\pi_i - \pi_j)$ is located in the $m$-th segment, e.g., $m = 5$ in Fig. 1. For segment $m$, the minimum errors of $g_l$ and $(\pi_i - \pi_j)$ are zero (i.e., at both ends of the segment), and the maximum errors of $g_l$ and $(\pi_i - \pi_j)$ are derived by the following method. Specifically, for gas passive pipeline $l$, let

$$g_l = \mathrm{a}_0(\pi_i - \pi_j) + \mathrm{b}_0.$$

Coefficient $\mathrm{a}_0 = (h(\mathrm{X}_l^{m+1}) - h(\mathrm{X}_l^m))/(\mathrm{X}_l^{m+1} - \mathrm{X}_l^m)$ is the gradient of segment $m$. Coefficient $\mathrm{b}_0$ is derived by jointly solving the following two equations and requiring a double root.

$$g_l = \mathrm{a}_0(\pi_i - \pi_j) + \mathrm{b}_0,$$

$$g_l^2 = \begin{cases} \mathrm{W}_l(\pi_i - \pi_j) & \text{if } \pi_i \ge \pi_j \\ -\mathrm{W}_l(\pi_i - \pi_j) & \text{if } \pi_i < \pi_j. \end{cases}$$

After determining $\mathrm{b}_0$, we can immediately know the maximum errors of $g_l$ and $(\pi_i - \pi_j)$. Here, we take Fig. 1 as an example, where $m = 5$. The coordinate of the red dot, (true $(\pi_i - \pi_j)$, true $g_l$), is the true values of $(\pi_i - \pi_j)$ and $g_l$ derived by original Weymouth equation (4j), and the coordinates of grey x mark and dark orange triangle, (estimated $(\pi_i - \pi_j)$, true $g_l$) and (true $(\pi_i - \pi_j)$, estimated $g_l$), are the PWL approximation of $g_l$ and $(\pi_i - \pi_j)$, respectively. Accordingly, we derive the maximum errors of $g_l$ and $(\pi_i - \pi_j)$, i.e., (estimated $(\pi_i - \pi_j)$ − true $(\pi_i - \pi_j)$) and (true $g_l$ − estimated $g_l$), respectively. By repeating the above procedures, the maximum errors of the PWL approxi-mation in other segments can be derived.

For an arbitrary segment that $(\pi_i - \pi_j)$ is located in, e.g., any one of the eight segments in Fig. 1, the minimum errors of $g_l$ and $(\pi_i - \pi_j)$ are also zero, and the maximum errors of $g_l$ and $(\pi_i - \pi_j)$ are the maximum values among all maximum errors of $g_l$ and $(\pi_i - \pi_j)$ in different segments, respectively.

For convenience, the compact form of the proposed bilevel LR model targeting IEGSs (4a)-(4c), (4e)-(4i), (4k)-(4o), (5), (6), denoted as the original model (**O-M**), is presented below.

$$\textbf{O-M}: \ \Gamma_{\text{O-M}} = \max_{\mathbf{x},\, \mathbf{y}_0,\, \mathbf{z}_0} \ \mathbf{c}^{\text{T}}\mathbf{y}_0 \tag{7a}$$

$$\text{s.t. } \mathbf{A}\mathbf{x} \le \mathbf{a}, \ \mathbf{x} \in \mathbb{R}^n, \tag{7b}$$

$$(\mathbf{y}_0,\ \mathbf{z}_0) \in \arg\min_{\mathbf{y},\, \mathbf{z}} \mathbf{c}^{\text{T}}\mathbf{y} \tag{7c}$$

$$\text{s.t. } \mathbf{E}\mathbf{y} + \mathbf{F}\mathbf{z} \le \mathbf{d} - \mathbf{D}\mathbf{x}, \ \mathbf{y} \in \mathbb{R}^q, \ \mathbf{z} \in \{0,\ 1\}^r. \tag{7d}$$

$\mathbf{x}$ denotes upper-level variables. $\mathbf{y}$ and $\mathbf{z}$ denote lower-level continuous and discrete variables, respectively. $\mathbf{A}$, $\mathbf{D}$, $\mathbf{E}$, and $\mathbf{F}$ are constant coefficient matrices. $\mathbf{a}$, $\mathbf{c}$, and $\mathbf{d}$ are constant vectors. For **O-M**, this paper assumes that there is at least one non-zero $\mathbf{x}$ satisfying (7b). Otherwise, intruders cannot launch any LR attack. This assumption is mild, as an IEGS usually contains more than one power/gas node that is connected with power/gas loads. Load redistribution manipulation among these power/gas loads is always feasible, leading to non-zero $\mathbf{x}$ satisfying (7b).

### *D. Model Extension*

This subsection reveals the extendibility of the proposed bilevel LR model **O-M** in the following aspects.

1) *Constraint violation*. The LR attack obtained by **O-M** and (3) may cause constraint violation and be suspected, as only the stealthiness of the attack is guaranteed. Adding (3d), (3e), and (8) to the upper level of **O-M** can effectively avoid this issue.

$$-\text{P}_l \le \sum_k \beta_{lk} \Big( \sum_{g \in \mathcal{P}_g \cup \mathcal{G}_g} \text{b}^{\text{G}}_{kg} \tilde{\text{p}}_g + \sum_{d \in \mathcal{P}_d} \text{b}^{\text{D}}_{kd} (\text{P}_d + \Delta p_d - \tilde{\text{s}}^{\text{p}}_d) \Big) \le \text{P}_l, \ \ \forall l \in \mathcal{P}_l, \tag{8a}$$

$$\text{G}^{\min}_n \le \tilde{\pi}_n + \Delta\pi_n \le \text{G}^{\max}_n, \ \ \forall n \in \mathcal{G}_n, \tag{8b}$$

$$\tilde{\pi}_{n(c)} + \Delta\pi_{n(c)} \le \alpha_c (\tilde{\pi}_{m(c)} + \Delta\pi_{m(c)}), \ \ \forall c \in \mathcal{G}_c, \tag{8c}$$

$$-\text{G}_l \le \tilde{\text{g}}_l + \Delta g_l \le \text{G}_l, \ \ \forall l \in \mathcal{G}_l, \tag{8d}$$

$$0 \le \tilde{\text{g}}_c + \Delta g_c \le \text{G}_c, \ \ \forall c \in \mathcal{G}_c, \tag{8e}$$

where $\tilde{\text{p}}_g$, $\tilde{\text{s}}^{\text{p}}_d$, $\tilde{\text{g}}_l$, $\tilde{\text{g}}_c$, and $\tilde{\pi}_n$ are the estimated measurements for $p_g$, $s^{\text{p}}_d$, $g_l$, $g_c$, and $\pi_n$, respectively. The nonconvex terms in (3e) can be linearized by the incremental PWL method [37], and the upper level of the extended **O-M** becomes an MILP.

2) *Attack budget*. In practice, intruders may consider attack budgets. A plausible method is to add constraints (1c), (3d), (3e), and (9) to the upper level of **O-M**.

$$\Delta p_d = 0 \Leftrightarrow \delta^{\text{p}}_d = 0, \ \ \forall d \in \mathcal{P}_d, \tag{9a}$$

$$\Delta p_l = 0 \Leftrightarrow \delta^{\text{p}}_l = 0, \ \ \forall l \in \mathcal{P}_l, \tag{9b}$$

$$\Delta g_d = 0 \Leftrightarrow \delta^{\text{g}}_d = 0, \ \ \forall d \in \mathcal{G}_d, \tag{9c}$$

$$\Delta \pi_n = 0 \Leftrightarrow \delta^{\text{g}}_n = 0, \ \ \forall n \in \mathcal{G}_n, \tag{9d}$$

$$\Delta g_l = 0 \Leftrightarrow \delta^{\text{g}}_l = 0, \ \ \forall l \in \mathcal{G}_l, \tag{9e}$$

$$\Delta g_c = 0 \Leftrightarrow \delta^{\text{g}}_c = 0, \ \ \forall c \in \mathcal{G}_c, \tag{9f}$$

$$\sum_{d \in \mathcal{P}_d} \delta^{\text{p}}_d + 2\sum_{l \in \mathcal{P}_l} \delta^{\text{p}}_l + \sum_{d \in \mathcal{G}_d} \delta^{\text{g}}_d + \sum_{n \in \mathcal{G}_n} \delta^{\text{g}}_n + \sum_{l \in \mathcal{G}_l} \delta^{\text{g}}_l + \sum_{c \in \mathcal{G}_c} \delta^{\text{g}}_c \le \text{B}, \tag{9g}$$

where $\delta^{\text{p}}_d$, $\delta^{\text{p}}_l$, $\delta^{\text{g}}_d$, $\delta^{\text{g}}_n$, $\delta^{\text{g}}_l$, and $\delta^{\text{g}}_c$ are binary auxiliary variables representing if the corresponding measurements are attacked or not (1: attacked, and 0 otherwise). Parameter B in (9g) refers to the attack budget. After linearizing (3e) (by the incremental PWL method [37]) and (equivalently) transforming (9a)-(9f) into mixed-integer linear constraints [11], the upper level of the extended **O-M** becomes an MILP.

3) *Power-to-gas (P2G) facilities*. IEGSs may install P2G facilities to reduce energy dissipation. Generally, the mathematical formulation of P2G facilities [26] is

$$g_{f^{\text{g}}} = \gamma_{f^{\text{p}}} p_{f^{\text{p}}}, \ \ \forall f^{\text{p}} \in \mathcal{P}_f, \ f^{\text{g}} \in \mathcal{G}_f,$$

where $p_{f^{\text{p}}}$ and $g_{f^{\text{g}}}$ denote the electricity transmitted to and the gas produced by P2G facilities, respectively. $\gamma_{f^{\text{p}}}$ is the conversion ratio of P2G facility $f$. $\mathcal{P}_f$ and $\mathcal{G}_f$ are the sets of P2G facilities connected to power and gas systems, respectively. Similar to power generators and gas wells, the amount of electricity transmitted to P2G facilities might be under the supervision of P2G control rooms and thus are easily detectable when false data are injected [11] (due to the direct communication between P2G control rooms and control centers). Considering that LR attacks are a kind of practical FDIA, intruders should not inject any false data into $p_{f^{\text{p}}}$ and $g_{f^{\text{g}}}$.

## III. Solution Method

The proposed **O-M** cannot be directly solved due to its bilevel structure. In this section, we design a modified R&D algorithm (based on the R&D algorithm [25]) to solve **O-M**.

### *A. Single-Level Reformulation*

According to [28], the relatively complete response assumption is critical for solving a bilevel MILP. For **O-M**, it requires that its lower-level problem is always feasible for any possible tuple ($\mathbf{x}$, $\mathbf{z}$). Any possible $\mathbf{x}$ and $\mathbf{z}$ refer to the $\mathbf{x}$ satisfying (7b) and the $\mathbf{z}$ satisfying (6), respectively. Under this assumption, a bilevel MILP can be equivalently reformulated as a single-level optimization problem, based on either the KKT conditions or strong duality theorem [28]. Unfortunately, **O-M** does not meet this requirement. A counterexample is the attack *col*(-0.6, 0.6) and the UC *col*(1, 1) in the illustrative example in Appendix A. In addition, the binary variables in (6) may also incur the violation of the relatively complete response assumption, rendering some attack strategies x in the upper of O-M infeasible. On the other hand, we find that some $\mathbf{z}$'s in **O-M** naturally satisfy this assumption. Thus, we divide all possible $\mathbf{z}$'s (in **O-M)** into two sets. Specifically, let set $\mathcal{Z}^{\mu}$ consist of all the possible $\mathbf{z}$'s in **O-M** that satisfy the relatively complete response assumption, i.e., $\mathcal{Z}^{\mu} = \{\mathbf{z}^{\mu*}_1, \mathbf{z}^{\mu*}_2, \cdots, \mathbf{z}^{\mu}_{|\mathcal{Z}^{\mu}|}{}^*\}$, and the rest of the possible $\mathbf{z}$'s in **O-M** constitute set $\mathcal{Z}^{\upsilon}$, i.e., $\mathcal{Z}^{\upsilon} = \{\mathbf{z}^{\upsilon*}_1, \mathbf{z}^{\upsilon*}_2, \cdots, \mathbf{z}^{\upsilon}_{|\mathcal{Z}^{\upsilon}|}{}^*\}$. Set $\mathcal{Z}$ is composed of all the possible $\mathbf{z}$'s in **O-M**, i.e., $\mathcal{Z} = \mathcal{Z}^{\mu} \cup \mathcal{Z}^{\upsilon}$, and is a finite set. $|\mathcal{Z}|$, $|\mathcal{Z}^{\mu}|$, and $|\mathcal{Z}^{\upsilon}|$ denote the number of elements in $\mathcal{Z}$, $\mathcal{Z}^{\mu}$, and $\mathcal{Z}^{\upsilon}$, respectively. Accordingly, we present a new KKT-based single-level reformulation for **O-M**, i.e., **KKT-R**, as follows.

$$\textbf{KKT-R}:\ \Gamma^*_{\text{K-M}} = \max_{\mathbf{x},\, \mathbf{y}_{|\mathcal{Z}|},\, \mathbf{z}_0,\, \boldsymbol{\mu}_{|\mathcal{Z}^\mu|},\, \boldsymbol{\upsilon}_{|\mathcal{Z}^\upsilon|},\, \mathbf{s}_{|\mathcal{Z}^\upsilon|}} \mathbf{c}^{\mathrm{T}}\mathbf{y}_0 \tag{10a}$$

$$\text{s.t. } \mathbf{Ax} \le \mathbf{a}, \tag{10b}$$

$$\mathbf{Dx} + \mathbf{Ey}_0 + \mathbf{Fz}_0 \le \mathbf{d}, \tag{10c}$$

$$\mathbf{c}^{\mathrm{T}}\mathbf{y}_0 \le \mathbf{c}^{\mathrm{T}}\mathbf{y}_i^\mu,\ \ 1 \le i \le |\mathcal{Z}^\mu|, \tag{10d}$$

$$(\mathbf{y}_i^\mu,\ \boldsymbol{\mu}_i) \in \mathcal{U}(\mathbf{x},\ \mathbf{z}_i^{\mu*}),\ \ 1 \le i \le |\mathcal{Z}^\mu|, \tag{10e}$$

$$\mathbf{c}^{\mathrm{T}}\mathbf{y}_0 \le \mathbf{c}^{\mathrm{T}}\mathbf{y}_j^\upsilon + \boldsymbol{\rho}^{\mathrm{T}}\mathbf{s}_j,\ \ 1 \le j \le |\mathcal{Z}^\upsilon|, \tag{10f}$$

$$(\mathbf{y}_j^\upsilon,\ \boldsymbol{\upsilon}_j,\ \mathbf{s}_j) \in \mathcal{V}(\mathbf{x},\ \mathbf{z}_j^{\upsilon*}),\ \ 1 \le j \le |\mathcal{Z}^\upsilon|, \tag{10g}$$

$$\mathbf{x} \in \mathbb{R}^n,\ \mathbf{y}_0 \in \mathbb{R}^q,\ \mathbf{z}_0 \in \{0,\ 1\}^r, \tag{10h}$$

$$\mathbf{y}_i^\mu \in \mathbb{R}^q,\ \boldsymbol{\mu}_i \in \mathbb{R}^p,\ \ 1 \le i \le |\mathcal{Z}^\mu|, \tag{10i}$$

$$\mathbf{y}_j^\upsilon \in \mathbb{R}^q,\ \boldsymbol{\upsilon}_j \in \mathbb{R}^p,\ \mathbf{s}_j \in \mathbb{R}_+^p,\ \ 1 \le j \le |\mathcal{Z}^\upsilon|. \tag{10j}$$

$\mathbf{y}_{|\mathcal{Z}|} = col(\mathbf{y}_0, \mathbf{y}_i^\mu, \mathbf{y}_j^\upsilon)$, $\boldsymbol{\mu}_{|\mathcal{Z}^\mu|} = col(\boldsymbol{\mu}_i)$, $\boldsymbol{\upsilon}_{|\mathcal{Z}^\upsilon|} = col(\boldsymbol{\upsilon}_j)$, $\mathbf{s}_{|\mathcal{Z}^\upsilon|} = col(\mathbf{s}_j)$, $1 \le i \le |\mathcal{Z}^\mu|$, $1 \le j \le |\mathcal{Z}^\upsilon|$. $\boldsymbol{\mu}_i$ is the dual variable for the lower-level problem of **O-M**, i.e., (7c), (7d), with fixed $\mathbf{z}_i^{\mu*}$. $\boldsymbol{\upsilon}_j$ is the dual variable for the relaxed lower-level problem of **O-M**, i.e., (11), with fixed $\mathbf{z}_j^{\upsilon*}$, where

$$\Gamma_r(\mathbf{x},\ \mathbf{z}_j^{\upsilon*}) = \min_{\mathbf{y}_j^\upsilon,\, \mathbf{s}_j} \mathbf{c}^{\mathrm{T}}\mathbf{y}_j^\upsilon + \boldsymbol{\rho}^{\mathrm{T}}\mathbf{s}_j \tag{11a}$$

$$\text{s.t. } \mathbf{Ey}_j^\upsilon - \mathbf{s}_j \le \mathbf{d} - \mathbf{Dx} - \mathbf{Fz}_j^{\upsilon*},\ \mathbf{y}_j^\upsilon \in \mathbb{R}^q,\ \mathbf{s}_j \in \mathbb{R}_+^p. \tag{11b}$$

$\mathbf{s}_j$ is the nonnegative slack variable. $\boldsymbol{\rho}$ is a (sufficiently) large positive penalty parameter. Sets

$$\mathcal{U}(\mathbf{x},\ \mathbf{z}_i^{\mu*}) = \begin{Bmatrix} \mathbf{Ey}_i^\mu \le \mathbf{d} - \mathbf{Dx} - \mathbf{Fz}_i^{\mu*}, \\ \mathbf{E}^{\mathrm{T}}\boldsymbol{\mu}_i + \mathbf{c} = \mathbf{0}, \\ \boldsymbol{\mu}_i \perp (\mathbf{Ey}_i^\mu + \mathbf{Dx} + \mathbf{Fz}_i^{\mu*} - \mathbf{d}), \\ \mathbf{0} \le \boldsymbol{\mu}_i, \end{Bmatrix} \tag{12a}$$

$$\mathcal{V}(\mathbf{x},\ \mathbf{z}_j^{\upsilon*}) = \begin{Bmatrix} \mathbf{Ey}_j^\upsilon - \mathbf{s}_j \le \mathbf{d} - \mathbf{Dx} - \mathbf{Fz}_j^{\upsilon*}, \\ \mathbf{E}^{\mathrm{T}}\boldsymbol{\upsilon}_j + \mathbf{c} = \mathbf{0}, \\ \boldsymbol{\upsilon}_j \perp (\mathbf{Ey}_j^\upsilon - \mathbf{s}_j + \mathbf{Dx} + \mathbf{Fz}_j^{\upsilon*} - \mathbf{d}), \\ (\boldsymbol{\rho} - \boldsymbol{\upsilon}_j) \perp \mathbf{s}_j, \\ \mathbf{0} \le \boldsymbol{\upsilon}_j \le \boldsymbol{\rho}. \end{Bmatrix} \tag{12b}$$

**KKT-R** is obtained by first enumerating all the possible $\mathbf{z}$'s in (7d) which satisfy the relatively complete response assumption (denoted as $\mathbf{z}_i^{\mu*}$, $\mathbf{z}_i^{\mu*} \in \mathcal{Z}^\mu$, where $\Gamma_f(\mathbf{z}_i^{\mu*}) = 0$). For these $\mathbf{z}_i^{\mu*}$, $\mathbf{z}_i^{\mu*} \in \mathcal{Z}^\mu$, we directly reformulate (7c), (7d) as (10d), (10e). Then, for the rest of the possible $\mathbf{z}$'s in (7d) which do not satisfy the relatively complete response assumption (denoted as $\mathbf{z}_j^{\upsilon*}$, $\mathbf{z}_j^{\upsilon*} \in \mathcal{Z}^\upsilon$, where $\Gamma_f(\mathbf{z}_j^{\upsilon*}) > 0$), we reformulate (7c), (7d) as (10f), (10g) by using the KKT-based reformulations of the relaxed lower-level problem (11). Consequently, **KKT-R** is always feasible for any possible $\mathbf{x}$. The following proposition states the relation between **O-M** and its reformulation **KKT-R**.

*Proposition 2*: If Assumption 1 holds and the optimal solution for **O-M** exists, there exists an optimal solution for **KKT-R** that is also optimal for **O-M**.

Please refer to Appendix B for the proof. Proposition 2 indicates an alternative for **O-M**. Namely, we can derive the optimal solution for **O-M** by solving **KKT-R**. Under Assumption 1, an optimal solution (or an $\epsilon$-optimal solution if the exact solution does not exist) for **O-M** is derived.

### *B. Modified R&D Algorithm*

The single-level reformulation **KKT-R** is computationally demanding due to a large number of elements in $\mathcal{Z}^\mu$ and $\mathcal{Z}^\upsilon$. On the other hand, it enables us to solve the bilevel **O-M** in a reformulation and decomposition manner [28]. Accordingly, we develop the modified R&D algorithm.

**Algorithm 1**: Modified R&D algorithm for solving **O-M**

**Step 1** Set $LB = -\infty$, $UB = +\infty$, $k = 0$, $k_\mu = 0$, $k_\upsilon = 0$, $\varepsilon > 0$.

**Step 2** Solve the master problem **MP** with given $\mathbf{z}_i^{\mu*}$, $1 \le i \le k_\mu$ and $\mathbf{z}_j^{\upsilon*}$, $1 \le j \le k_\upsilon$.

$$\textbf{MP}:\ \bar{\Gamma} = \max_{\mathbf{x},\, \mathbf{y}_k,\, \mathbf{z}_0,\, \boldsymbol{\mu}_k,\, \boldsymbol{\upsilon}_k,\, \mathbf{s}_k} \mathbf{c}^{\mathrm{T}}\mathbf{y}_0 \tag{13a}$$

$$\text{s.t. } \mathbf{Ax} \le \mathbf{a}, \tag{13b}$$

$$\mathbf{Dx} + \mathbf{Ey}_0 + \mathbf{Fz}_0 \le \mathbf{d}, \tag{13c}$$

$$\mathbf{c}^{\mathrm{T}}\mathbf{y}_0 \le \mathbf{c}^{\mathrm{T}}\mathbf{y}_i^\mu,\ \ 1 \le i \le k_\mu, \tag{13d}$$

$$(\mathbf{y}_i^\mu,\ \boldsymbol{\mu}_i) \in \mathcal{U}(\mathbf{x},\ \mathbf{z}_i^{\mu*}),\ \ 1 \le i \le k_\mu, \tag{13e}$$

$$\mathbf{c}^{\mathrm{T}}\mathbf{y}_0 \le \mathbf{c}^{\mathrm{T}}\mathbf{y}_j^\upsilon + \boldsymbol{\rho}^{\mathrm{T}}\mathbf{s}_j,\ \ 1 \le j \le k_\upsilon, \tag{13f}$$

$$(\mathbf{y}_j^\upsilon,\ \boldsymbol{\upsilon}_j,\ \mathbf{s}_j) \in \mathcal{V}(\mathbf{x},\ \mathbf{z}_j^{\upsilon*}),\ \ 1 \le j \le k_\upsilon, \tag{13g}$$

$$\mathbf{x} \in \mathbb{R}^n,\ \mathbf{y}_0 \in \mathbb{R}^q,\ \mathbf{z}_0 \in \{0,\ 1\}^r, \tag{13h}$$

$$\mathbf{y}_i^\mu \in \mathbb{R}^q,\ \boldsymbol{\mu}_i \in \mathbb{R}^p,\ \ 1 \le i \le k_\mu, \tag{13i}$$

$$\mathbf{y}_j^\upsilon \in \mathbb{R}^q,\ \boldsymbol{\upsilon}_j \in \mathbb{R}^p,\ \mathbf{s}_j \in \mathbb{R}_+^p,\ \ 1 \le j \le k_\upsilon. \tag{13j}$$

Obtain the optimal LR attack decision $\mathbf{x}^*$ and the optimal objective value $\bar{\Gamma}$. Set $UB = \bar{\Gamma}$.

**Step 3** If $UB$ - $LB \le \varepsilon$, return $\mathbf{x}^*$ and $\bar{\Gamma}$, and then terminate; otherwise, go to Step 4.

**Step 4** Solve the subproblem 1 (**SP1**) with given $\mathbf{x}^*$.

$$\textbf{SP1}:\ \underline{\Gamma}(\mathbf{x}^*) = \min_{\mathbf{y},\, \mathbf{z}} \mathbf{c}^{\mathrm{T}}\mathbf{y} \tag{14a}$$

$$\text{s.t. } \mathbf{Ey} + \mathbf{Fz} \le \mathbf{d} - \mathbf{Dx}^*, \tag{14b}$$

$$\mathbf{y} \in \mathbb{R}^q,\ \mathbf{z} \in \{0,\ 1\}^r. \tag{14c}$$

Obtain the optimal solution $\mathbf{z}^*$ and the optimal objective value $\underline{\Gamma}(\mathbf{x}^*)$. Set $LB = \max(LB, \underline{\Gamma}(\mathbf{x}^*))$.

**Step 5** Solve the subproblem 2 (**SP2**) with given $\mathbf{z}^*$.

$$\textbf{SP2}:\ \Gamma_f(\mathbf{z}^*) = \max_{\mathbf{x}} \min_{\mathbf{y},\, \mathbf{s}} \mathbf{1}^{\mathrm{T}}\mathbf{s} \tag{15a}$$

$$\text{s.t. } \mathbf{Ax} \le \mathbf{a}, \tag{15b}$$

$$\mathbf{Dx} + \mathbf{Ey} - \mathbf{s} \le \mathbf{d} - \mathbf{Fz}^*, \tag{15c}$$

$$\mathbf{x} \in \mathbb{R}^n,\ \mathbf{y} \in \mathbb{R}^q,\ \mathbf{s} \in \mathbb{R}_+^p. \tag{15d}$$

Obtain the optimal objective value $\Gamma_f(\mathbf{z}^*)$. If $\Gamma_f(\mathbf{z}^*) = 0$, set $\mathbf{z}^\mu_{k_\mu+1}{}^* = \mathbf{z}^*$, and update $k_\mu = k_\mu+1$; otherwise, set $\mathbf{z}^\upsilon_{k_\upsilon+1}{}^* = \mathbf{z}^*$, and update $k_\upsilon = k_\upsilon+1$. Update $k = k + 1$, and go to Step 2.

†: Constraints (13d), (13e), and (13i) are not incorporated in **MP** if $k_\mu = 0$, and constraints (13f), (13g), and (13j) are not incorporated in **MP** if $k_\upsilon = 0$.

$\varepsilon$ is the convergence threshold. In **MP**, $\mathbf{y}_k = col(\mathbf{y}_0, \mathbf{y}_i^\mu, \mathbf{y}_j^\upsilon)$, $1 \le i \le k_\mu$, $1 \le j \le k_\upsilon$, $\boldsymbol{\mu}_k = col(\boldsymbol{\mu}_i)$, $1 \le i \le k_\mu$, $\boldsymbol{\upsilon}_k = col(\boldsymbol{\upsilon}_j)$, and $\mathbf{s}_k = col(\mathbf{s}_j)$, $1 \le j \le k_\upsilon$. As is shown in Fig. 2, the modified R&D algorithm is implemented iteratively in a master-subproblem framework. First, we initialize the parameters presented in Step 1. At each iteration, we first derive an optimal $\mathbf{x}^*$ by solving **MP**, the *relaxation* of **KKT-R**. Then, we evaluate if the convergence condition in Step 3 is satisfied. If yes, we derive the optimal LR attack decision $\mathbf{x}^*$ and terminate; otherwise, we solve **SP1**, the lower-level problem of **O-M** with fixed $\mathbf{x}^*$, and obtain the optimal $\mathbf{z}^*$. Next, we use **SP2** to check the feasibility of this $\mathbf{z}^*$ for all possible $\mathbf{x}$. Note that it is the proposed **SP2** that achieves dividing $\mathbf{z}^*$ into two sets. For any $\mathbf{x}$, the objective function of **SP2** is to minimize the summation of the slack variables in lower-level constraints. The objective value must be equal to zero if no constraint is violated and larger than zero if at least one constraint is violated. That is why $\mathbf{z}^* \in \mathcal{Z}^{\mu}$ if $\Gamma_f(\mathbf{z}^*) = 0$ and $\mathbf{z}^* \in \mathcal{Z}^{\upsilon}$ if $\Gamma_f(\mathbf{z}^*) > 0$. Specifically, if $\Gamma_f(\mathbf{z}^*) = 0$, indicating $\mathbf{z}^* \in \mathcal{Z}^{\mu}$, cuts (13d), (13e), and (13i) are added to **MP**, and we turn to the next iteration; if $\Gamma_f(\mathbf{z}^*) > 0$, indicating $\mathbf{z}^* \in \mathcal{Z}^{\upsilon}$, cuts (13f), (13g), and (13j) are added to **MP**, and we turn to the next iteration. The convergence of this algorithm is analyzed in the following proposition.

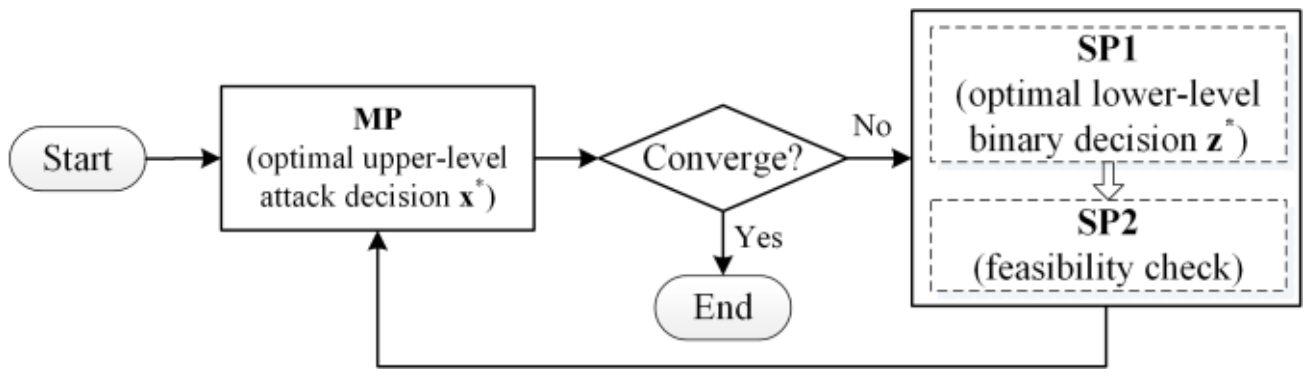

Fig. 2. Framework of the modified R&D algorithm.

*Proposition 3*: If Assumption 1 holds and the optimal solution for **O-M** exists, the modified R&D algorithm converges to the optimum of **O-M** within $O(|\mathcal{Z}|)$ iterations.

The proof of Proposition 3 follows the fact that the set $\mathcal{Z}$ is a finite set and any repeated $\mathbf{z}_i^{\mu *}$ or $\mathbf{z}_j^{\upsilon *}$ in $\mathcal{Z}$ will lead to $UB = LB$ (Proposition 8 in [28]). If the optimal solution for **O-M** does not exist, this algorithm converges to the $\epsilon$-optimum of **O-M** under Assumption 1. Note that Algorithm 1 also applies to solving the extended models in Section II.D.

TABLE IV
COMPARISON BETWEEN MODEL (16) AND THE KKT-BASED REFORMULATION FOR **SP2**

| | No. of variables | Difference | No. of constraints | Difference |
|---|---|---|---|---|
| Model (16) | $m + n + p$ | - | $3m + n + 2p + q$ | - |
| KKT-M | $n + 2p + q$ | $p + q - m$ | $m + 6p + q$ | $4p - 2m - n$ |

Before the implementation of Algorithm 1, we must fix two issues. First, **MP** is nonconvex due to the bilinear terms in sets $\mathcal{U}$ and $\mathcal{V}$, i.e., the third constraint in (12a) and the third and fourth constraints in (12b). We adopt the Big-M method [28] to linearize them. Second, the bilevel **SP2** is computationally intractable. This paper proposes a computationally tractable and efficient reformulation for **SP2**. Based on both the strong duality theorem and KKT conditions, we reformulate **SP2** as

$$\max_{\mathbf{x},\,\boldsymbol{\pi},\,\boldsymbol{\kappa}} \ (\mathbf{F}\mathbf{z}^* - \mathbf{d})^{\mathrm{T}}\boldsymbol{\pi} + \mathbf{a}^{\mathrm{T}}\boldsymbol{\kappa} \tag{16a}$$

$$\text{s.t. } \mathbf{A}\mathbf{x} \le \mathbf{a}, \tag{16b}$$

$$\mathbf{E}^{\mathrm{T}}\boldsymbol{\pi} = \mathbf{0}, \tag{16c}$$

$$\mathbf{0} \le \boldsymbol{\pi} \le \mathbf{1}, \tag{16d}$$

$$\mathbf{D}^{\mathrm{T}}\boldsymbol{\pi} - \mathbf{A}^{\mathrm{T}}\boldsymbol{\kappa} = \mathbf{0}, \tag{16e}$$

$$\boldsymbol{\kappa} \perp (\mathbf{a} - \mathbf{A}\mathbf{x}), \tag{16f}$$

$$\boldsymbol{\kappa} \ge \mathbf{0}, \tag{16g}$$

$$\mathbf{x} \in \mathbb{R}^n,\ \boldsymbol{\pi} \in \mathbb{R}^p,\ \boldsymbol{\kappa} \in \mathbb{R}^m. \tag{16h}$$

Please refer to Appendix C for more details. The Big-M method is adopted to linearize (16f), converting (16) into an MILP, which can be solved by off-the-shelf solvers, such as Gurobi. As an alternative, the (pure) KKT-based method also applies to reformulating **SP2** as a (single-level) MILP (short for KKT-M). We compare these two reformulations for **SP2** and present the results in Table IV. The definitions of $m$, $n$, $p$, and $q$ are given in (15d) and (16h), where $p > m > n$. Compared with the KKT-based reformulation (KKT-M), model (16) has fewer variables and constraints (the differences in numbers are shown in Table IV) and is expected to have better computational efficiency.

## IV. CASE STUDY

This section validates the effectiveness of the proposed bi-level LR attack model **O-M** and solution method using an integrated 6-node-electricity-7-node-gas system (IEGS-6-7) and an integrated 39-node-electricity-20-node-gas system (IEGS-39-20). Reference [38] provides detailed system data and topology. Tests are coded using Julia 1.6.2 with JuMP toolbox and conducted on a laptop with an Intel(R) i5-8265U CPU @ 1.60 GHz and a 16 GB memory. MILPs are solved by Gurobi 9.1.2.

### A. *Integrated 6-Node-Electricity-7-Node-Gas System*

The IEGS-6-7 contains one coal-fired generator ($G_1$), two gas-fired generators ($G_2$ and $G_3$), three power loads ($PL_1$, $PL_2$, and $PL_3$), two gas wells ($GW_1$ and $GW_2$), and three gas loads ($GL_1$, $GL_2$, and $GL_3$). The output costs satisfy $G_3 > G_2 > G_1$ and $GW_2 > GW_1$, and the gas load shedding costs are set to be higher than the power load shedding costs. The number of piecewise line segments K is set to 8. The penalty parameter $\boldsymbol{\rho}$ in **MP** is set to 3e1. The Big-M values in **MP** and **SP2** are set to 1e2. The convergence threshold ε in Step 3 of Algorithm 1 is set to 1e-4, and the solver's tolerance is set to 1e-5.

TABLE V
LR ATTACKS UNDER DIFFERENT CONDITIONS

| Scenario | | A1 | A2 | A3 | A4 | A5 | A6 |
|---|---|---|---|---|---|---|---|
| $\tau^p$ and $\tau^g$ | | 0% | 30% | 50% | 50% | 50% | 50% |
| Load | Power (MW) | 352 | 352 | 352 | **370** | 252 | 52 |
| | Gas (Sm$^3$/h) | 5324 | 5324 | 5324 | **4800** | 4324 | 2324 |
| Objective value ($) | | 7099.6 | 7151.5 | 7592.4 | 8567.4 | 4987.4 | 1538.6 |
| LR attack loss | | 0.0% | 0.7% | 6.9% | **18.2%** | **0.0%** | **0.0%** |
| Generator output (MW) | $G_1$ | 220.0 | 216.8 | 212.3 | 209.8 | 220.0 | **52.0** |
| | $G_2$ | 100.0 | 92.3 | 83.4 | 90.0 | 32.0 | **0.0** |
| | $G_3$ | 32.0 | 42.9 | 50.0 | 50.0 | **0.0** | **0.0** |
| Gas well output (Sm$^3$/h) | $GW_1$ | 4924.0 | 5163.5 | 5162.0 | 4900.0 | 0.0 | 0.0 |
| | $GW_2$ | 6000.0 | 6000.0 | 6000.0 | 6000.0 | 5604.0 | 2324.0 |
| Load shedding | Power (MW) | 0.0 | 0.0 | **6.3** | **20.2** | 0.0 | 0.0 |
| | Gas (Sm$^3$/h) | 0.0 | 0.0 | 0.0 | 0.0 | 0.0 | 0.0 |

1) *LR* attacks *on IEGSs*. This part tests the impact of LR attacks on IEGSs under different load forecasting errors, i.e., $\tau^p$ and $\tau^g$ (reflecting the impacts of load forecasting accuracy on LR attacks), and load levels. Power and gas transmission limits are fixed to 125 MW and 6000 Sm$^3$/h, respectively. Table V shows the results. Scenario A1 is the normal IEGS operation state without LR attacks and is selected as the base case for calculating LR attack losses of Scenarios A2 and A3. Test results show that LR attacks can lead to i) a slight increase in system costs by affecting generator and gas well outputs *without load shedding* (Scenario A2), or ii) a significant increase in system costs by incurring load-shedding (Scenario A3). Next, we fix $\tau^p$ and $\tau^g$ to 50% and test the impact of LR attacks under different load levels (Scenarios A4-A6). Compared with Scenario A3, the LR attack on heavy load Scenario A4 induces much higher costs. When we reduce loads to the level of Scenario A5, none of LR attacks cause economic losses, i.e., cost increase. We continue to decrease loads to the level of Scenario A6, where the aggregated minimum generator output (70 MW) is higher than the aggregated power loads (52 MW). Benefitted from UC variables, **O-M** can still provide LR attack strategies, while the bilevel LR attack model (4) cannot due to the violation of minimum generation outputs. Overall, this test demonstrates that i) LR attacks on the IEGSs with poor load forecasting accuracy and heavy loads are liable to induce higher economic losses; ii) incorporating UC variables enables the proposed LR attack model **O-M** to cope with more IEGS operation states and is considered significant.

TABLE VI
SOLUTION ACCURACY

| Gas passive pipeline model | Obj. value ($) | Generator output (MW) | | | Gas well output ($Sm^3/h$) | | Load shedding | |
|---|---|---|---|---|---|---|---|---|
| | | $G_1$ | $G_2$ | $G_3$ | $GW_1$ | $GW_2$ | Power (MW) | Gas ($Sm^3/h$) |
| Inc. PWL | 7592.4 | 212.3 | 83.4 | 50.0 | 5162.0 | 6000.0 | 6.3 | 0.0 |
| Ori. Wey. | 7592.4 | 212.3 | 83.5 | 50.0 | 5162.0 | 6000.0 | 6.3 | 0.0 |

In addition, we test the impact of the incremental PWL method on solution accuracy, and the testing results are presented in Table VI. Scenario A3 is chosen as the base case for using both the incremental PWL method and original Weymouth equations. Obj., Inc., Ori., and Wey. are short for objective, incremental, original, and Weymouth, respectively. To the best of our knowledge, no existing method can solve a bilevel nonconvex optimization problem with a nonconvex lower level. The results of using Ori. Wey are derived by fixing upper-level LR attack values (derived by Inc. PWL), transforming the bilevel MILP into a (single- level) MILP. According to Table VI, the incremental PWL method does not affect the objective value, generator outputs, gas well outputs, or load shedding at all. Only the values of gas nodal pressures have some slight differences. Please refer to [38] for the detailed test data. This test numerically validates our error analysis in Section II.C and the effectiveness of the incremental PWL method for simplifying Weymouth equations.

TABLE VII
LR ATTACKS UNDER DIFFERENT TRANSMISSION LIMITS

| Scenario | | B1 | B2 | B3 | B4 |
|---|---|---|---|---|---|
| Transmission limit | $P_l$ (MW) | 125 | 110 | 150 | 200 |
| | $G_l$ and $G_c$ ($Sm^3/h$) | 7000 | 5500 | 8000 | 10000 |
| Objective value under normal state ($) | | 6571.6 | 6707.1 | 6571.6 | 6571.6 |
| Objective value under LR attack ($) | | 6653.0 | 7370.7 | 6571.6 | 6571.6 |
| LD attack loss | | **1.2%** | **9.9%** | **0.0%** | **0.0%** |

2) *Influence of transmission limits on LR attacks*. This part investigates the impact of LR attacks on IEGSs under different transmission limits, and testing results are shown in Table VII. Transmission limits refer to the $P_l$ in (4e), $G_l$ in (4l), and $G_c$ in (4m). The aggregated power and gas loads are fixed to 302 MW and 6324 $Sm^3/h$, respectively, and $\tau^p$ and $\tau^g$ are fixed to 50%. Scenario B1 is selected as the base case. When we shrink transmission limits to the level of Scenario B2, the LR attack leads to larger economic losses. When transmission limits are relaxed to 150 MW and 8000 $Sm^3/h$ (Scenario B3), none of LR attacks can incur increases in system costs. Interestingly, we find that LR attacks on lossless IEGSs cannot affect system costs when all transmission limits are inactive. A lossless IEGS is an IEGS that can transmit electricity and gas without any loss, e.g., the proposed **O-M**. Transmission limits are considered inactive if i) constraints (4e), (4l), and (4m) are inactive, and ii) nodal pressures are sufficient for gas flow to reach transmission limits. Physically, if an IEGS is able to transmit an arbitrary amount of electricity and gas to any power bus and gas node in this IEGS without any loss, the distribution of loads cannot affect system costs. This is the reason why LR attacks become invalid for Scenarios B3, B4, A5, and A6. Thus, the LR attacks derived by **O-M** are also called *congestion-related cyberattacks*. This case indicates the IEGSs with heavy loads and limited transmission capability are more vulnerable to LR attacks.

TABLE VIII
COMPARISON BETWEEN R&D AND MODIFIED R&D ALGORITHMS

| Load | Power (MW) | 352 | | 252 | |
|---|---|---|---|---|---|
| | Gas ($Sm^3/h$) | 5324 | | 6324 | |
| Solution algorithm | | U-R&D | M-R&D | U-R&D | M-R&D |
| Objective value ($) | | 7513.8 | **7592.4** | 6190.3 | **6395.0** |
| No. of iterations | | 5 | 6 | 5 | 4 |
| Injected false data | $\Delta PL_1$ (MW) | 35.2 | -35.2 | 0.0 | 0.0 |
| | $\Delta PL_2$ (MW) | -70.4 | 70.4 | -100.8 | -100.8 |
| | $\Delta PL_3$ (MW) | 35.2 | -35.2 | 100.8 | 100.8 |
| | $\Delta GL_1$ ($Sm^3/h$) | -391.7 | 887.3 | -834.0 | 2108.0 |
| | $\Delta GL_2$ ($Sm^3/h$) | 29.0 | -1250 | 1388.0 | -3162.0 |
| | $\Delta GL_3$ ($Sm^3/h$) | 362.7 | 362.7 | -554.0 | 1054.0 |

3) *Importance of the relatively complete response assumption*. This part is to i) compare the proposed modified R&D algorithm (M-R&D) with the R&D algorithm without verifying the relatively complete response assumption (U-R&D) and ii) evaluate the impact of this assumption on solving the proposed LR attack model **O-M**, and Table VIII shows the results. $\Delta PL_i$ and $\Delta GL_i$, $i \in \{1, 2, 3\}$, refer to the injected false data into power and gas loads, respectively. Power and gas transmission limits are fixed to 125 MW and 6000 $Sm^3/h$, respectively. $\tau^p$ and $\tau^g$ are set to 50% and 100% for different load levels (i.e., 352 MW and 5324 $Sm^3/h$ and 252 MW and 6324 $Sm^3/h$), respectively. As is shown in Table VIII, U-R&D fails to obtain the optimal solutions for both load levels due to the *suboptimal LR attacks on gas loads*. Theoretically, if the relatively complete response assumption does not hold, U-R&D can only find a suboptimal LR attack. In contrast, the proposed M-R&D derives (global) optimal solutions for both load levels. In addition, we observe that U-R&D and M-R&D converge within a similar number of iterations for both load levels, indicating that their computation times are very close (all within 100 seconds). In fact, M-R&D is more complicated than U-R&D, as the former addresses feasibility issues. This test exhibits that the proposed SP2 does not incur heavy computational burdens for M-R&D, validating its effectiveness and efficiency.

TABLE IX
IMPACT OF MODEL EXTENSION ON LR ATTACKS

| Scenario | | Constraints | | Attack budgets | |
|---|---|---|---|---|---|
| Model | | **O-M** | Model extension 1) | **O-M** | Model extension 2) |
| Objective value under LR attack ($) | | **7513.8** | 7124.4 | **8567.4** | 8406.3 |
| Injected false data | $\Delta PL_1$ (MW) | 35.2 | -3.1 | -37.0 | -37.0 |
| | $\Delta PL_2$ (MW) | -70.4 | 70.4 | 74.0 | 37.0 |
| | $\Delta PL_3$ (MW) | 35.2 | -67.3 | -37.0 | 0.0 |
| | $\Delta GL_1$ ($Sm^3/h$) | 358.3 | 605.7 | 800.0 | 800.0 |
| | $\Delta GL_2$ ($Sm^3/h$) | -470.9 | 162.0 | -500.0 | -800.0 |
| | $\Delta GL_3$ ($Sm^3/h$) | 112.6 | 443.7 | -300.0 | 0.0 |

4) *Model extension cases*. This part tests the impact of model extensions (in Section II.D) on the consequence of LR attacks. Test results are presented in Table IX. Model extensions 1) and 2) refer to the bilevel LR attack models introduced in Section II.D 1) and 2), respectively. The gross power and gas loads are set to 352 MW and 5324 $Sm^3/h$ and 370 MW and 4800 $Sm^3/h$, respectively. Power and gas transmission limits are set to 125 MW and 7000 $Sm^3/h$ and 125 MW and 6000 $Sm^3/h$, respectively. $\tau^p$ and $\tau^g$ are set to 50%. According to Table IX, the losses of IEGS-6-7 caused by LR attacks decrease significantly when intruders incorporate all the constraints in model extension 1). On the other hand, this attack does not violate any operation constraints and is considered more stealthy. In practice, intruders may consider balancing the impact and the stealthiness of LR attacks by adding partial stealthy constraints. For model extension 2), we assume that intruders can only compromise two power loads (two out of three) and two gas loads (two out of three) due to limited attack capability. Unsurprisingly, the LR attack with attack budgets causes lower economic losses than that without attack budgets. Overall, this case presents the extendibility of the proposed bilevel LR attack model and the value of the extensions for addressing more sophisticated cases.

TABLE X
LR ATTACKS UNDER DIFFERENT CONDITIONS

| Scenario | | C1 | C2 | C3 | C4 |
|---|---|---|---|---|---|
| $\tau^p$ and $\tau^g$ | | 0% | 50% | 50% | 50% |
| Load | Power (MW) | 1650 | 1650 | 1750 | 1650 |
| | Gas ($Sm^3/h$) | 1200 | 1200 | 1400 | 1200 |
| Transmission limit | | 100%‡ | 100% | 100% | 90% |
| Objective value ($) | | 43847.2 | 50559.0 | 58367.0 | 55067.9 |
| LR attack loss | | 0% | 15.3% | 18.8% | 25.6% |
| No. of iterations | | 1 | 2 | 3 | 3 |
| Computation time (s) | | 29 | 102 | 597 | 608 |

‡: Please refer to [38] for the detailed (100%) transmission limit data.

### B. *Integrated 39-Node-Electricity-20-Node-Gas System*

This subsection tests the proposed bilevel LR attack model **O-M** and the solution method on the medium-sized IEGS-39-20. This IEGS has seven coal-fired generators, three gas-fired generators, nineteen power loads, two gas wells, and nine gas loads. The number of piecewise line segments K is set to 4. The penalty parameter $\boldsymbol{\rho}$ in **MP** is set to 1e3. The Big-M values in **MP** and **SP2** are set to 1e3. The convergence threshold $\varepsilon$ in Step 3 of Algorithm 1 is set to 1e-3, and the solver's tolerance is set to 1e-4. Table X shows the testing results. Computation time refers to the execution time of the entire program. Scenario C1, the normal IEGS operation state without LR attacks, is selected as the base case. Similar to the results in the last subsection, the impact of LR attacks on this medium-sized IEGS is affected by load forecasting errors (Scenario C2), load levels (Scenario C3), and transmission limits (Scenario C4). Moreover, all scenarios are solved within 608 seconds and a small number of iterations, which is reasonable for practical system dispatch (based on [39]). This test validates the scalability of the proposed bilevel LR attack model **O-M** and the solution method (i.e., the modified R&D algorithm).

## V. CONCLUSION

This paper investigates the paradigm and stealthy conditions of LR attacks on IEGSs. In order to identify the most severe LR attack (from an economic perspective), this paper proposes a bilevel MILP and develops a modified R&D algorithm to solve the proposed bilevel model. Testing results indicate that: i) LR attacks on the IEGSs with poor load forecasting accuracy and heavy loads are prone to induce higher economic losses; ii) incorporating UC variables enables the proposed LR attack model **O-M** to cope with more operation states and is considered significant; iii) LR attacks on lossless IEGSs are intrinsically congestion-related attacks and cannot affect system costs when all transmission limits are inactive; iv) the modified R&D algorithm is effective in addressing the bilevel MILPs when the relatively complete response assumption is not satisfied. Note that IEGS operators can also use the proposed bilevel LR attack model **O-M** to evaluate and identify IEGS cyber-vulnerability when facing LR attacks in an offline manner. However, this work does not consider i) cyberattacks on IEGSs with incomplete network information, ii) multi-period LR attack models with energy storage facilities and line packs, iii) the cascading effect of LR attacks on IEGSs, or iv) detection and mitigation strategies to detect these LR attacks and reduce their impact. These will be our future works.

## APPENDIX

### A. *An Illustrative Example*

We use a 2-bus power system (shown in Fig. 3) to illustrate the infeasibility issue. Power loads $PL_1$ and $PL_2$ are 3.5 MW and 6.5 MW, respectively. The output ranges of $G_1$ and $G_2$ are [5, 10] and [2, 5] MW, respectively. The thermal limit of PLine is 2 MW. As a big penalty, we set the load shedding costs to 300 $/MWh and 400 $/MWh for $PL_1$ and $PL_2$, respectively. Scenario 1 in Table XI presents the economic dispatch decision under the normal operation state (without any LR attack). $LS_1$ and $LS_2$ denote the load shedding in $PL_1$ and $PL_2$, respectively.

Assume that intruders launch an LR attack by setting $\widehat{PL}_1 =$ 2.5 MW and $\widehat{PL}_2 = 7.5$ MW (Scenario 2), where $\widehat{PL}_1$ and $\widehat{PL}_2$ are falsified power loads on buses 1 and 2, respectively. For this attack, there is no feasible dispatch decision, as it violates the minimum output limit of $G_1$. Due to the lower-level infeasibility, this attack is not a feasible solution for the bilevel model (4) and thus is ignored naturally by intruders. However, based on the definition in Section II.B, it is, undoubtedly, a possible LR attack. Under this attack, the operator has to shut down $G_1$, yielding a cost of $2,110 (Scenario 3), which compared to the worst LR attack derived by (5) (Scenario 4), is more costly.

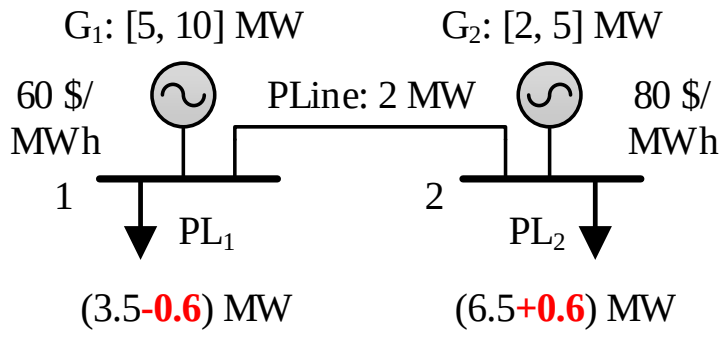

Fig. 3. An illustrative example.

TABLE XI
ECONOMIC DISPATCH DECISIONS UNDER DIFFERENT SCENARIOS

| LR attack scenario | System Cost ($) | $G_1$ (MW) | $G_2$ (MW) | $LS_1$ (MW) | $LS_2$ (MW) | PLine (MW) |
|---|---|---|---|---|---|---|
| Scenario 1 | **690** | 5.5 | 4.5 | 0 | 0 | 2 |
| Scenario 2 | - | - | - | - | - | - |
| Scenario 3 | **2110** | 0 | 5 | 2.9 | 2.1 | 0 |
| Scenario 4 | **700** | 5 | 5 | 0 | 0 | 2 |

## B. Proof of Proposition 2

*Proof*: First, we prove (B.1)-(B.7) is a relaxation of **KKT-R**.

$$\max_{\mathbf{x}'', \mathbf{y}''_{|\mathcal{Z}|}, \mathbf{z}''_0, \boldsymbol{\upsilon}''_{|\mathcal{Z}|}, \mathbf{s}''_{|\mathcal{Z}|}} \mathbf{c}^{\mathrm{T}}\mathbf{y}''_0 \tag{B.1}$$

$$\text{s.t. } \mathbf{A}\mathbf{x}'' \le \mathbf{a}, \tag{B.2}$$

$$\mathbf{D}\mathbf{x}'' + \mathbf{E}\mathbf{y}''_0 + \mathbf{F}\mathbf{z}''_0 \le \mathbf{d}, \tag{B.3}$$

$$\mathbf{c}^{\mathrm{T}}\mathbf{y}''_0 \le \mathbf{c}^{\mathrm{T}}\mathbf{y}''_j + \boldsymbol{\rho}^{\mathrm{T}}\mathbf{s}''_j, \quad 1 \le j \le |\mathcal{Z}|, \tag{B.4}$$

$$(\mathbf{y}''_j, \boldsymbol{\upsilon}''_j, \mathbf{s}''_j) \in \mathcal{V}(\mathbf{x}'', \mathbf{z}''^{*}_j), \quad 1 \le j \le |\mathcal{Z}|, \tag{B.5}$$

$$\mathbf{x}'' \in \mathbb{R}^n, \mathbf{y}''_0 \in \mathbb{R}^q, \mathbf{z}''_0 \in \{0, 1\}^r, \tag{B.6}$$

$$\mathbf{y}''_j \in \mathbb{R}^q, \boldsymbol{\upsilon}''_j \in \mathbb{R}^p, \mathbf{s}''_j \in \mathbb{R}^p_+, \quad 1 \le j \le |\mathcal{Z}|. \tag{B.7}$$

$\mathbf{y}''_{|\mathcal{Z}|} = col(\mathbf{y}''_0, \mathbf{y}''_j)$, $\boldsymbol{\upsilon}''_{|\mathcal{Z}|} = col(\boldsymbol{\upsilon}''_j)$, $\mathbf{s}''_{|\mathcal{Z}|} = col(\mathbf{s}''_j)$, $1 \le j \le |\mathcal{Z}|$. $\{\mathbf{z}''^{*}_1, \mathbf{z}''^{*}_2, \cdots, \mathbf{z}''^{*}_{|\mathcal{Z}|}\} = \{\mathbf{z}^{\mu*}_1, \mathbf{z}^{\mu*}_2, \cdots, \mathbf{z}^{\mu}_{|\mathcal{Z}^{\mu}|}{}^{*}, \mathbf{z}^{\upsilon*}_1, \cdots, \mathbf{z}^{\upsilon}_{|\mathcal{Z}^{\upsilon}|}{}^{*}\}$. Model (B.1)-(B.7) is derived by replacing (10d), (10e), (10i) in **KKT-R** with (10f), (10g), (10j). Assume that $(\mathbf{x}, \mathbf{y}_{|\mathcal{Z}|}, \mathbf{z}_0, \boldsymbol{\mu}_{|\mathcal{Z}^{\mu}|}, \boldsymbol{\upsilon}_{|\mathcal{Z}^{\upsilon}|}, \mathbf{s}_{|\mathcal{Z}^{\upsilon}|})$ is a feasible solution for **KKT-R**. Let $(\mathbf{x}'', \mathbf{y}''_{|\mathcal{Z}|}, \mathbf{z}''_0, \boldsymbol{\upsilon}''_{|\mathcal{Z}|}, \mathbf{s}''_{|\mathcal{Z}|}) = (\mathbf{x}, \mathbf{y}_{|\mathcal{Z}|}, \mathbf{z}_0, (\boldsymbol{\mu}_{|\mathcal{Z}^{\mu}|}, \boldsymbol{\upsilon}_{|\mathcal{Z}^{\upsilon}|}), (\mathbf{0}_{|\mathcal{Z}^{\mu}|}, \mathbf{s}_{|\mathcal{Z}^{\upsilon}|}))$, where $\mathbf{0}_{|\mathcal{Z}^{\mu}|}$ is a $(p * |\mathcal{Z}^{\mu}|) * 1$ all zero vector. Evidently, $(\mathbf{x}'', \mathbf{y}''_{|\mathcal{Z}|}, \mathbf{z}''_0, \boldsymbol{\upsilon}''_{|\mathcal{Z}|}, \mathbf{s}''_{|\mathcal{Z}|})$ is a feasible solution for (B.1)-(B.7), as model (B.1)-(B.7) degenerates to **KKT-R** if we set the first $p * |\mathcal{Z}^{\mu}|$ elements in $\mathbf{s}''_{|\mathcal{Z}|}$ to zero. Since any feasible solution for **KKT-R** is also feasible for (B.1)-(B.7), (B.1)-(B.7) is a relaxation of **KKT-R**.

Similarly, we can prove that **KKT-R** is a relaxation of **O-M**. In addition, it is easy to prove that (B.1)-(B.7) is an equivalent single-level reformulation of (B.8)-(B.12).

$$\max_{\mathbf{x}'', \mathbf{y}''_0, \mathbf{z}''_0} \mathbf{c}^{\mathrm{T}}\mathbf{y}''_0 \tag{B.8}$$

$$\text{s.t. } \mathbf{A}\mathbf{x}'' \le \mathbf{a}, \ \mathbf{x}'' \in \mathbb{R}^n, \tag{B.9}$$

$$(\mathbf{y}''_0, \mathbf{z}''_0) \in \arg\min_{\mathbf{y}'', \mathbf{z}'', \mathbf{s}''} \mathbf{c}^{\mathrm{T}}\mathbf{y}'' + \boldsymbol{\rho}^{\mathrm{T}}\mathbf{s}'' \tag{B.10}$$

$$\text{s.t. } \mathbf{E}\mathbf{y}'' + \mathbf{F}\mathbf{z}'' - \mathbf{s}'' \le \mathbf{d} - \mathbf{D}\mathbf{x}'', \tag{B.11}$$

$$\mathbf{y}'' \in \mathbb{R}^q, \ \mathbf{z}'' \in \{0, 1\}^r, \ \mathbf{s}'' \in \mathbb{R}^p_+. \tag{B.12}$$

According to Proposition 5 in [28], (i.e., there exists an optimal solution for (B.8)-(B.12) that is also feasible and optimal for **O-M**,) we immediately conclude that there exists an optimal solution for **KKT-R** that is also optimal for **O-M**. This completes the proof. ■

## C. Single-Level Reformulation for SP2

Based on the strong duality theorem, **SP2** is equivalent to

$$\max_{\mathbf{x}, \boldsymbol{\pi}} \boldsymbol{\pi}^{\mathrm{T}}(\mathbf{D}\mathbf{x} + \mathbf{F}\mathbf{z}^{*} - \mathbf{d}) \tag{C.1}$$

$$\text{s.t. } \mathbf{A}\mathbf{x} \le \mathbf{a}, \tag{C.2}$$

$$\mathbf{E}^{\mathrm{T}}\boldsymbol{\pi} = \mathbf{0}, \tag{C.3}$$

$$\mathbf{0} \le \boldsymbol{\pi} \le \mathbf{1}, \tag{C.4}$$

$$\mathbf{x} \in \mathbb{R}^n, \boldsymbol{\pi} \in \mathbb{R}^p. \tag{C.5}$$

We re-write (C.1)-(C.5) as

$$\max_{\boldsymbol{\pi} \in \Pi} \boldsymbol{\pi}^{\mathrm{T}}(\mathbf{F}\mathbf{z}^{*} - \mathbf{d}) + \max_{\mathbf{x} \in \mathcal{X}} \boldsymbol{\pi}^{\mathrm{T}}\mathbf{D}\mathbf{x}, \tag{C.6}$$

where $\Pi := \{\boldsymbol{\pi} \in \mathbb{R}^p \mid \text{(C.3)-(C.4)}\}$ and $\mathcal{X} := \{\mathbf{x} \in \mathbb{R}^n | \text{(C.2)}\}$. Since **SP2** is feasible for any $\mathbf{x}$ satisfying (7b), model (C.6) has a finite objective value. We employ the KKT conditions to reformulate the second term in (C.6), $\max_{\mathbf{x} \in \mathcal{X}} \boldsymbol{\pi}^{\mathrm{T}}\mathbf{D}\mathbf{x}$, as

$$\max_{\mathbf{x}, \boldsymbol{\kappa}} \boldsymbol{\pi}^{\mathrm{T}}\mathbf{D}\mathbf{x} \tag{C.7}$$

$$\text{s.t. } \mathbf{A}\mathbf{x} \le \mathbf{a}, \tag{C.8}$$

$$\boldsymbol{\pi}^{\mathrm{T}}\mathbf{D} - \boldsymbol{\kappa}^{\mathrm{T}}\mathbf{A} = \mathbf{0}, \tag{C.9}$$

$$\boldsymbol{\kappa}^{\mathrm{T}}(\mathbf{a} - \mathbf{A}\mathbf{x}) = \mathbf{0}, \tag{C.10}$$

$$\boldsymbol{\kappa} \ge \mathbf{0}, \tag{C.11}$$

$$\mathbf{x} \in \mathbb{R}^n, \boldsymbol{\kappa} \in \mathbb{R}^m. \tag{C.12}$$

Multiplying (C.9) by $\mathbf{x}$, we have

$$\boldsymbol{\pi}^{\mathrm{T}}\mathbf{D}\mathbf{x} - \boldsymbol{\kappa}^{\mathrm{T}}\mathbf{A}\mathbf{x} = \mathbf{0}. \tag{C.13}$$

According to (C.10) and (C.13), the bilinear objective (C.7), $\boldsymbol{\pi}^{\mathrm{T}}\mathbf{D}\mathbf{x}$, can be replaced by $\boldsymbol{\kappa}^{\mathrm{T}}\mathbf{a}$, and we derive model (16).